\documentclass[12pt,reqno]{amsart}
\usepackage[utf8]{inputenc}
\usepackage{amsmath}
\usepackage{amsfonts}
\usepackage{amssymb}
\usepackage{amsthm}
\usepackage{hyperref}
\usepackage{multicol, comment}
\usepackage{pictexwd, dcpic}
\usepackage{fancyhdr}
\usepackage{mathrsfs}
\usepackage{latexsym,mathtools}
\usepackage{bm, enumerate}
\usepackage{graphicx}
\usepackage{color,soul}
\usepackage{xcolor}
\usepackage{dsfont}
\usepackage{todonotes}
\usepackage{tikz-cd}

\makeatletter
\@namedef{subjclassname@2020}{%
	\textup{2020} Mathematics Subject Classification}
\makeatother

\newcommand{\beq}{\begin{eqnarray}}
	\newcommand{\eeq}{\end{eqnarray}}
\newcommand{\beqn}{\begin{eqnarray*}}
	\newcommand{\eeqn}{\end{eqnarray*}}

\reversemarginpar

\numberwithin{equation}{section}

\newtheorem{theorem}{\bf Theorem}[section]
\newtheorem{proposition}[theorem]{\bf Proposition}
\newtheorem{corollary}[theorem]{\bf Corollary}
\newtheorem{lemma}[theorem]{\bf Lemma}

\theoremstyle{remark}
\newtheorem{definition}[theorem]{\bf Definition}
\newtheorem{example}[theorem]{\bf Example}
\newtheorem{Non-example}[theorem]{\bf Non-Example}
\newtheorem{remark}[theorem]{\bf Remark}

\newcommand{\bydef}{\stackrel{\text{def}}{=}}

\newcommand\T{\mathbb{T}}
\newcommand{\D}{\mathbb{D}}
\newcommand{\N}{\mathbb{N}}
\newcommand{\C}{\mathbb{C}}
\newcommand{\Z}{\mathbb{Z}}

\theoremstyle{definition}

\theoremstyle{remark}

\hypersetup{
	colorlinks,
	citecolor=blue,
	linkcolor=black,
	urlcolor=blue}

\begin{document}
\title[Numerical Radius Interpolation]{Interpolation, Factorization, and Approximation in the Numerical Radius Unit Ball}

\author[Ahmed]{Hamidul Ahmed}
\address{Department of Mathematics, Indian Institute of Technology Bombay, Powai, Mumbai, 400076, India}
\email{hamidulahmedslk@gmail.com, 22d0785@iitb.ac.in}

\author[Bhowmik]{Mainak Bhowmik}
\address{Department of Mathematics, Indian Institute of Technology Bombay, Powai, Mumbai, 400076, India}
\email{mainakbhowmik943@gmail.com, mainak@math.iitb.ac.in}

\author [Das]{B. Krishna Das}
\address{Department of Mathematics, Indian Institute of Technology Bombay, Powai, Mumbai, 400076, India}
\email{bata436@gmail.com, dasb@math.iitb.ac.in}

\author[Sau]{Haripada Sau}
\address{Department of Mathematics, Indian Institute of Science Education and Research, Dr.\ Homi Bhabha Road, Pashan, Pune, Maharashtra 411008, India.}
\email{haripadasau215@gmail.com, hsau@iiserpune.ac.in}



\begin{abstract}
We initiate a systematic study of Nevanlinna--Pick interpolation and function theory for matrix-valued holomorphic functions on the unit disk taking values in the numerical radius unit ball---a class we term \emph{$\nu$-Schur functions}.

In contrast to the operator-norm setting, where commutant lifting reduces interpolation to the positivity of a single Pick matrix, we show that a natural one-parameter family of Pick-type matrices arising from Ando’s characterization is necessary but not sufficient for solvability.
We overcome this obstruction by establishing two necessary and sufficient solvability criteria. The first is constructive, resting on a factorization of a $\nu$-Schur function $F = PGQ$ into a classical Schur function $G$ and bounded holomorphic matrix functions $P, Q$ coupled by a boundary condition. This yields an accompanying Choi--Li-type dilation representation. The second is dual and geometric, characterizing solvability via a Nehari-type distance from a canonical model-space interpolant to a subspace determined by the interpolation data.

Finally, we introduce \emph{$\nu$-inner functions}---$\nu$-Schur functions whose non-tangential boundary values lie almost everywhere in the set of extremal numerical contractions---and prove that every non-extremal problem admits a $\nu$-inner solution. As a consequence, we obtain a numerical-radius Carath\'eodory approximation theorem: $\nu$-inner functions are dense in the $\nu$-Schur class under the compact-open topology.
\end{abstract}

\maketitle

\tableofcontents

\section{Introduction}

The theory of interpolation by holomorphic functions stands at the confluence of complex analysis, operator theory, and system theory. The foundational Nevanlinna--Pick interpolation problem asks whether there exists a Schur function -- that is, a holomorphic function $f: \mathbb{D} \to \mathbb{D}$ on the open unit disk $\mathbb{D} = \{z \in \mathbb{C} : |z| < 1\}$ -- satisfying prescribed interpolation conditions $f(\lambda_j) = w_j$ for distinct nodes $\lambda_1, \dots, \lambda_N \in \mathbb{D}$ and targets $w_1, \dots, w_N \in \mathbb{D}$. Pick's classic 1915 theorem established that such an interpolant exists if and only if the Pick matrix 
\[
\left[ \frac{1 - w_j \bar{w}_k}{1 - \lambda_j \bar{\lambda}_k} \right]_{j,k=1}^N \succeq 0
\]
is positive semi-definite \cite{Pick} while Nevanlinna in 1919 paratmetrized all the solutions in \cite{Nevanlinna}. The Pick interpolation problem has inspired a vibrant and far-reaching field of research over the past century, serving as a fundamental catalyst for modern operator theory and function theory. Half a century after Pick's work, Sarason \cite{Sarason1967} re-conceptualized Nevanlinna--Pick interpolation as an operator-theoretic statement on compressed shift operators acting on the model space $H^2 \ominus \theta H^2$, a shift-operator perspective that culminated in the celebrated commutant lifting theorem. Parallel efforts to extend interpolation beyond the open unit disk led to groundbreaking advances on multiply connected domains and Riemann surfaces, pioneered by Abrahamse and Douglas \cite{Abrahamse1979, Abrahamse-Douglas1976} through the theory of bundle shifts, and systematically expanded by Ball and his collaborators -- most notably with Helton \cite{Ball-Helton1983}, Gohberg and Rodman \cite{Ball-Gohberg-Rodman1990}, Clancey \cite{Ball-Clancey1992}, and Trent \cite{Ball-Trent1998} -- to matrix and operator-valued interpolation across complex manifolds, Riemann surfaces, and multivariable domains. Subsequently, Agler \cite{Agler-OpThAdvAppl} extended interpolation theory to the bidisk and polydisk, replacing standard Schur functions with what is now known as the Schur--Agler class. Ball and Bolotnikov extended this perspective for domains defined by matrix polynomials \cite{Ball-Bolotnikov}. Motivated by these novel ideas, Dritschel and McCullough \cite{DM} developed an abstract test-function machinery for operator-valued interpolation across general operator systems, demonstrating that positivity with respect to appropriate families of test functions serves as the overarching paradigm for operator-norm interpolation.

In this article, we introduce and solve a new paradigm: the {\it numerical radius interpolation problem}. Given distinct nodes $\boldsymbol{\Lambda} = (\lambda_1, \ldots, \lambda_N) \in \mathbb{D}^N$ and target $d\times d$ complex matrices $\boldsymbol{W} = (W_1, \ldots, W_N) \in M_d^N$, we seek a holomorphic function $f \in H^\infty(\mathbb{D}, M_d)$ satisfying
\begin{equation} \label{eq:IP_intro}
f(\lambda_j) = W_j \quad \text{for } 1 \le j \le N, \quad \text{and} \quad \|f\|_{\nu, \infty} \bydef \sup_{z \in \mathbb{D}} \nu(f(z)) \le 1,
\end{equation}
where $\nu(T) \bydef \sup_{\|x\|=1} |\langle Tx, x \rangle|$ denotes the numerical radius of a matrix $T \in M_d$.

A primary conceptual motivation for our investigation stems from the failure of Sarason-type commutant lifting techniques when operator norm contractivity $\|T\| \le 1$ is replaced by numerical radius contractivity $\nu(T) \le 1$. Recall that an operator $T$ is a contraction if and only if $I - T^*T \succeq 0$. The positivity of $I - T^*T$ yields the necessary and sufficient Pick matrix condition for classical operator-valued interpolation. On the other hand, an operator $T$ is a {\it numerical contraction}, i.e., $\nu(T) \le 1$ if and only if
\[
I_d - \operatorname{Re}(\zeta T) \succeq 0 \quad \text{for every } \zeta \in \mathbb{T} =  \{z \in \mathbb{C} : |z|=1\}.
\]
While $I - T^* T \succeq 0$ seamlessly induces positive semi-definite kernel matrices on Hilbert space, the continuous family of conditions $I - \operatorname{Re}(\zeta T) \succeq 0$ provides only a necessary condition for the solvability of a numerical radius interpolation problem. Indeed, while Theorem \ref{thm:Nec-cond} shows that the conditions
\begin{align}\label{Eq:Pick-family}
    P(\theta) \bydef \left[ \frac{2 I_d - e^{i\theta} W_j - e^{-i\theta} W_k^*}{1 - \lambda_j \bar{\lambda}_k} \right]_{j,k=1}^N \succeq 0 \quad \mbox{for all }\theta \in [0, 2\pi),
\end{align} are necessary for the $N$-point numerical radius interpolation problem $\lambda_j\mapsto W_j$ to be solvable, Example \ref{exp:roadblock} shows that these conditions are not sufficient, in general. This failure underscores a key structural obstruction: a compressed numerical contraction on a co-invariant subspace does not, in general, lift to a $\nu$-contractive multiplication operator on the Hardy space. See Remark~\ref{rm:CLT} for a precise and elaborated discussion to this effect.

Beyond this operator-theoretic dichotomy lies a second, equally formidable roadblock of a function-theoretic and geometric nature. The Schur algorithm recursively reduces the degree of a matrix Nevanlinna-Pick interpolation problem by composing candidate functions with matrix Möbius transformations. The efficacy of the Schur algorithm hinges on the transitivity of the matrix-norm unit ball $\{ A \in M_d(\mathbb{C}) : \|A\| < 1 \}$ under its holomorphic automorphism group, enabling any interior target matrix to be mapped to the origin while preserving the norm constraint. For the numerical radius norm, however, the geometry of the domain is far more rigid: in the forthcoming paper \cite{ABDS}, the authors establish that for $d \ge 2$, the biholomorphic automorphism group of the numerical-radius unit ball 
\begin{equation}\label{Nu-Ball}
\overline{\Omega}_\nu \bydef \{ T \in M_d(\mathbb{C}) : \nu(T) \le 1 \}.
\end{equation}
does not act transitively on it. In fact, every such automorphism is linear and hence, the orbit of the zero matrix is singleton.
The inhomogeneity of the numerical-radius unit ball impedes the construction of degree-reducing linear fractional transformations, rendering standard Schur-type inductive schemes ineffective for numerical radius interpolation and necessitating non-traditional approach.

Interpolation under non-standard constraints has a distinguished history. In particular, Bercovici, Foias, and Tannenbaum \cite{Bercovici} addressed the spectral radius interpolation problem, replacing the standard operator norm bound $\|f\|_\infty \le 1$ with the spectral radius constraint $r(f(z)) \le 1$ through spectral commutant lifting techniques. In this article, not only do we introduce a novel interpolation theory based on the numerical radius, but also we provide two completely different solvability criteria.

The first solvability criterion (Theorem~\ref{thm:interpolation}) yields a direct, constructive characterization rooted in a structural factorization of holomorphic functions taking values in $\overline{\Omega}_\nu$, which we refer to as \emph{$\nu$-Schur functions}. In Section~\ref{sec:structure}, we prove that a function $f \in H^\infty(\mathbb{D}, M_d)$ is a $\nu$-Schur function if and only if it admits the factorization
\begin{equation*}
f = \Phi_1 G \Phi_2,
\end{equation*}
where $G \in H^\infty(\mathbb{D}, M_d)$ is a classical matrix-valued Schur function ($\|G\|_\infty \bydef \sup_{z\in\mathbb{D}}\|G(z)\| \le 1$) and $\Phi_1, \Phi_2 \in H^\infty(\mathbb{D}, M_d)$ are holomorphic matrix functions whose non-tangential boundary values satisfy the quadratic envelope condition
\begin{equation}\label{GenIneq}
\Phi_1(\zeta) \Phi_1(\zeta)^* + \Phi_2(\zeta)^* \Phi_2(\zeta) \preceq 2I_d \quad \text{for a.e. } \zeta \in \mathbb{T}.
\end{equation}
Equivalently, $f$ takes the form
\begin{align}\label{DilationYesIntro}
    f = V^*\begin{bmatrix}
        0 & 2G\\
        0 & 0
    \end{bmatrix}V,
\end{align}
where $G$ is a Schur function and $V:\mathbb{D}\to\mathcal{B}(\mathbb{C}^d,\mathbb{C}^{2d})$ is a contractive map formed from holomorphic and conjugate-holomorphic matrix functions. Consequently, the non-tangential limit $V^\flat(\zeta)\bydef\lim_{r\to 1-}V(r\zeta)$ exists almost everywhere on $\mathbb{T}$. Although we expect that $V^\flat$ takes isometric values almost everywhere on $\mathbb{T}$, we establish this under the strict bound $\|f\|_{\nu,\infty}<1$. A central analytical tool behind these representations is the spectral factorization of matrix-valued \emph{spectral densities} on $\mathbb{T}$ into outer factors.

Notably, the representation \eqref{DilationYesIntro} provides a function-theoretic counterpart to a classical theorem of Choi and Li~\cite{Choi}, which characterizes numerical contractions $T$ on a Hilbert space $\mathcal{H}$ ($\nu(T)\leq 1$) as
$$
T = V^*\begin{bmatrix}
    0 & 2I_\mathcal{H}\\
    0 & 0
\end{bmatrix}V
$$
for some isometry $V:\mathcal{H}\to\mathcal{H}\oplus\mathcal{H}$.

In contrast to this constructive factorization approach, our second solvability criterion (Theorem~\ref{thm:distance-criterion}) adopts a dual, geometric perspective through a Nehari-type distance formula. Passing to the predual space $L_{\nu^*}^1$, endowed with a norm $\nu^*$ dual to the numerical radius, we express the reciprocal of the minimal $\nu$-supremum norm 
$$
 \gamma_{\boldsymbol{\Lambda}, \boldsymbol{W}} \bydef \inf \{ \|F\|_{\nu,\infty}: \ F \in H^\infty(\D, M_d) \text{ and } F(\lambda_j)=W_j \text{ for } 1\leq j \leq N \}.
$$ as the $L_{\nu^*}^1$-distance from a normalized interpolant $\psi_{\boldsymbol{\Lambda},\boldsymbol{W}} / \|\psi_{\boldsymbol{\Lambda},\boldsymbol{W}}\|_2^2$ to a subspace $\mathcal{M}_{\boldsymbol{\Lambda},\boldsymbol{W}}$ both canonically determined by the interpolation data.

No interpolation theory is complete without a companion approximation theory it illuminates. Indeed, this paper is as much a foundational study of holomorphic functions taking values in the numerical radius unit ball as it is an investigation of the numerical radius interpolation problem itself. In classical function theory, a celebrated paradigm originating with Carath\'eodory asserts that every Schur function on $\mathbb{D}$ can be approximated locally uniformly by finite Blaschke products \cite{Caratheodory}. The Blaschke products are unimodular on the unit circle and hence, are rational inner functions. A proof, which seamlessly extends itself to matrix-valued Schur functions, of Carath\'eodory's theorem can be done via the existence of matrix-valued rational inner solution, for a solvable Nevanlinna--Pick problem on several domains such as the disk, the bidisk, the quotients of polydisk \cite{AM-Crelle, BK-JFA, DKS, Sarason1967}. 

 We establish numerical-radius counterparts of these landmark results. Central to this development is the notion of \emph{$\nu$-inner functions}: $\nu$-Schur functions whose non-tangential boundary values belong almost everywhere to the set of \emph{extremal numerical contractions}---a distinguished subset of the topological boundary of $\overline{\Omega}_\nu$ introduced and studied by Dritschel and Woerdeman~\cite{DW}. Theorem~\ref{thm:strict} proves that every \textit{non-extremal} numerical radius interpolation problem, i.e., a problem for which $\gamma_{\boldsymbol{\Lambda}, \boldsymbol{W}} <1$,  admits a $\nu$-inner solution. As a consequence, Theorem~\ref{thm:main} yields a numerical-radius Carath\'eodory-type approximation theorem: the collection of $\nu$-inner functions is dense in the $\nu$-Schur class with respect to the compact-open (Fr\'echet) topology on $\mathcal{O}(\mathbb{D}, M_d)$, the space of all $M_d$-valued holomorphic functions on $\mathbb{D}$. Ultimately, these results demonstrate that $\nu$-inner functions play the same foundational and extremal role in numerical radius function theory as Blaschke products and classical inner functions do in operator-norm analysis.

\section{A natural necessary condition that fails to be sufficient} \label{sec:roadblock}

Throughout the paper we say that $T\in M_d$ is a {\it numerical contraction} if $$\nu(T)=\sup_{x \in \C^d, \ \|x\|=1}|\langle Tx,x\rangle|\leq 1.$$ 
Moreover, $T$ is called a strict numerical contraction whenever $\nu(T)<1$. 

A matrix $T$ in $M_d$ is a contraction if and only if $I_d-T^*T\succeq0$. The following is a direct analogue of this for numerical contractions. See \cite{Ando, Gustafson-Rao1997} or \cite[Lemma~1.2]{DW} for a proof.

\begin{proposition} \label{prop:NRC-chrac}
A matrix $T$ in $M_d$ is a numerical contraction if and only if 
\begin{align}\label{Delta} 
 I_d - \operatorname{Re} (\zeta T) \succeq 0 \quad \mbox{for every $\zeta \in \mathbb{T}$}.
\end{align}
\end{proposition}

Although the numerical radius interpolation problem was formulated conceptually in the Introduction (see the discussion surrounding \eqref{eq:IP_intro}), it constitutes the central object of study throughout this manuscript. To fix terminology and notation for the node and target tuples used in our technical analysis, we formally restate the problem below.

\begin{definition}\label{def:IP}
Given an $N$-tuple of pairwise distinct nodes $\boldsymbol{\Lambda} = (\lambda_1, \ldots, \lambda_N) \in \mathbb{D}^N$ and an $N$-tuple of target matrices $\boldsymbol{W} = (W_1, \ldots, W_N) \in M_d^N$, the {\it Numerical Radius Interpolation Problem} asks for the existence of a holomorphic function $f \in H^\infty(\mathbb{D}, M_d)$ satisfying $\|f\|_{\nu, \infty} \le 1$ and
\begin{align} \label{eq:IP_def}
f(\lambda_j) = W_j, \quad j=1, \dots, N.
\end{align}
For brevity, we shall refer to this problem as a {\it $\nu$-interpolation problem}, or the {\it $\nu$-interpolation problem $\boldsymbol{\Lambda} \mapsto \boldsymbol{W}$} when the underlying data needs to be emphasized.
\end{definition}
As the title of the section indicates, we shall find a necessary condition for a $\nu$-interpolation problem to be solvable that fails to be sufficient. Proposition \ref{prop:NRC-chrac} provides the necessary condition, which, in turn, has two equivalent features that are derived by a different characterization of numerical contractions as follows. We shall see in later sections that it is this characterization that provides the fundamental operator-theoretic leverage required to overcome the obstacle demonstrated in this section.
\begin{proposition}[Theorem~1 in Ando \cite{Ando}]\label{prop:ando}
A matrix $T$ is a numerical contraction if and only if there exists a self-adjoint matrix $Z$ such that
\begin{equation}\label{eq:ando-block}
  \begin{bmatrix}
    2I_d - Z & T^*\\
    T & Z
  \end{bmatrix} \succeq 0.
\end{equation}
\end{proposition}
We shall now present the promised necessary conditions that are not sufficient.
\begin{theorem} \label{thm:Nec-cond}
Let ${\bf \Lambda}=(\lambda_1, \ldots, \lambda_N)\mapsto \boldsymbol{W}=(W_1, \ldots, W_N)$ be a $\nu$-interpolation problem. 
Set
\[
S\bydef \left[\frac{I_d}{1-\lambda_i\overline{\lambda}_j}\right]_{i,j=1}^{N},
\quad
D \bydef \operatorname{diag}(W_1,\ldots,W_N).
\]
If the $\nu$-interpolation problem 
$\bf\Lambda \mapsto \boldsymbol{W}
$ is solvable, then
\begin{align}\label{Eq:Pick-family}
    P(\theta) = \left[ \frac{2 I_d - e^{i\theta} W_j - e^{-i\theta} W_k^*}{1 - \lambda_j \overline{\lambda}_k} \right]_{j,k=1}^N \succeq 0 \quad \mbox{for all }\theta \in [0, 2\pi).
\end{align} 
Moreover, the following are equivalent:
\begin{enumerate}
\item[\textup{(i)}] $P(\theta)\succeq0$ for every $\theta\in[0,2\pi)$;
\item[\textup{(ii)}] $\nu(S^{-1/2}DS^{1/2})\leq1$;
\item[\textup{(iii)}] there exists a self-adjoint block matrix
$Z=[Z_{ij}]_{i,j=1}^{N}\in M_{Nd}$ such that
\[
\begin{bmatrix}
S-Z_{\mathrm{s}} & DS\\
SD^* & S+Z_{\mathrm{s}}
\end{bmatrix}\succeq0,
\quad
Z_{\mathrm{s}}
\bydef \left[
\frac{Z_{ij}}{1-\lambda_i\overline{\lambda}_j}
\right]_{i,j=1}^{N}.
\]
\end{enumerate}
\end{theorem}

\begin{proof} 
Suppose $f \in H^\infty(\mathbb{D}, M_d)$ solves the interpolation problem. So, $\nu(f(z)) \le 1$ for every $z \in \mathbb{D}$. Then for each $\theta \in [0, 2\pi)$, the function $G_\theta(z) = I_d - e^{i\theta} f(z)$ is in the Herglotz class on $\mathbb{D}$ i.e., $\operatorname{Re}(G_\theta(z)) \succeq 0$ for $z \in \mathbb{D}$. Therefore, the matrix-valued kernel
\begin{align*}
   K_\theta(z, w) \bydef \frac{G_\theta(z) + G_\theta(w)^*}{1 - z \bar{w}} \ \ \text{for } z, w \in \mathbb{D},
\end{align*}
is positive semi-definite on $\mathbb{D} \times \mathbb{D}$ for any $\theta \in [0, 2\pi)$; see \cite{Agler-OpThAdvAppl}. For any $\theta \in [0, 2\pi)$, the positive semi-definiteness of $K_\theta$ yields:
\begin{align}\label{eq:positivity-Pick}
\sum_{j, k=1}^N \left\langle \frac{G_\theta(\lambda_j) + G_\theta(\lambda_k)^*}{1 - \lambda_j \overline{\lambda}_k} u_k, u_j \right\rangle = \sum_{j, k=1}^N \left\langle \frac{2I_d - e^{i\theta} W_j - e^{-i\theta} W_k^*}{1 - \lambda_j \overline{\lambda}_k} u_k, u_j \right\rangle \ge 0.
\end{align}
Hence, $P(\theta) \succeq 0$ for all $\theta \in [0, 2\pi)$.

We now prove the equivalence. Since the Szeg\"{o} kernel functions at distinct points are linearly independent, its Gram matrix is positive definite. Thus $S$ is invertible. Set
\[
T \bydef S^{-1/2}DS^{1/2}.
\]
Then using the identity $P(\theta)=2S-e^{i\theta}DS-e^{-i\theta}SD^*$, we note that
\[
S^{-1/2}P(\theta)S^{-1/2}
=2I_{Nd}-e^{i\theta}T-e^{-i\theta}T^*
=2\bigl(I_{Nd}-\operatorname{Re}(e^{i\theta}T)\bigr).
\]
This proves the equivalence of \textup{(i)} and \textup{(ii)}.

Suppose that \textup{(ii)} holds. Since $\nu(T^*)=\nu(T)\leq1$, Proposition~\ref{prop:ando} provides a matrix $0\preceq L\preceq2I_{Nd}$ such that
\[
\begin{bmatrix}
2I_{Nd}-L & T\\
T^* & L
\end{bmatrix}\succeq0.
\]
Put
\[
K \bydef S^{1/2}(L-I_{Nd})S^{1/2},
\ \text{ and }
Z_{ij}=(1-\lambda_i\overline{\lambda}_j)K_{ij}.
\]
Since $K=K^*$, the block matrix $Z=[Z_{ij}]$ is self-adjoint,
and $Z_{\mathrm{s}}=K$. Also by this choice,
\[S-Z_{\mathrm{s}}=S-K=S^{1/2}(2I_{Nd}-L)S^{1/2}\ \text{ and } \ S+Z_{\mathrm{s}}=S^{1/2}L S^{1/2}.\]
Consequently,
\[
\begin{bmatrix}
S-Z_{\mathrm{s}} & DS\\
SD^* & S+Z_{\mathrm{s}}
\end{bmatrix}
=
\begin{bmatrix}
S^{1/2}&0\\
0&S^{1/2}
\end{bmatrix}
\begin{bmatrix}
2I_{Nd}-L&T\\
T^*&L
\end{bmatrix}
\begin{bmatrix}
S^{1/2}&0\\
0&S^{1/2}
\end{bmatrix}
\succeq0.
\]
Here we have used $S^{1/2}TS^{1/2}=DS$.
This proves \textup{(iii)}.

Finally, if \textup{(iii)} holds, then a straightforward matrix multiplication reveals that 
\[
P(\theta)
=
\begin{bmatrix}
I_{Nd} &
-e^{-i\theta}I_{Nd}
\end{bmatrix}
\begin{bmatrix}
S-Z_{\mathrm{s}} & DS\\
SD^* & S+Z_{\mathrm{s}}
\end{bmatrix}
\begin{bmatrix}
I_{Nd}\\
-e^{i\theta}I_{Nd}
\end{bmatrix}
\succeq0
\quad \text{for } \theta\in[0,2\pi).
\]
Thus \textup{(iii)} implies \textup{(i)}, completing the proof.
\end{proof}

We now demonstrate by a concrete example that the necessary condition obtained in Theorem \ref{thm:Nec-cond} is not sufficient, in general.

\begin{example}\label{exp:roadblock}
Consider the interpolation nodes and target matrices
\[
\lambda_1=0,\quad \lambda_2=\frac12,
\quad
W_1=\begin{bmatrix}0&1\\0&0\end{bmatrix},
\quad
W_2=\begin{bmatrix}0&0\\1&0\end{bmatrix}.
\]
Clearly, $\nu(W_1)=\nu(W_2)=1/2$. We show that the necessary
positivity condition holds, but the $\nu$-interpolation problem is not solvable.
Fix $\theta\in [0,2\pi)$, and write $\omega=e^{i\theta}\in\mathbb T$. Then a simple calculation reveals that 
\[
 P(\theta)
=
\left[
\frac{2I_2-\omega W_j-\overline{\omega}W_k^*}
     {1-\lambda_j\overline{\lambda_k}}
\right]_{j,k=1}^{2}
=
\begin{bmatrix}
A&B\\
B^*&C
\end{bmatrix},
\]
where
\[
A=
\begin{bmatrix}
2&-\omega\\
-\overline{\omega}&2
\end{bmatrix},
\quad
B=
\begin{bmatrix}
2&-(\omega+\overline{\omega})\\
0&2
\end{bmatrix},\
\text{ and }\
C=\frac43
\begin{bmatrix}
2&-\overline{\omega}\\
-\omega&2
\end{bmatrix}.
\]
Since $A$ has eigenvalues $1$ and $3$, it is positive definite. Moreover, using 
\[
A^{-1}
=\frac13
\begin{bmatrix}
2&\omega\\
\overline{\omega}&2
\end{bmatrix},
\]
a direct computation gives
\[
B^*A^{-1}B=C.
\]
Thus, the Schur complement is zero, and consequently
$P(\theta)\succeq0$. Since $\theta$ was arbitrary, the necessary positivity condition is met. Suppose, however, that there exists a holomorphic function
$F:\mathbb D\to M_2(\mathbb C)$ satisfying
\[
\sup_{z\in\mathbb D}\nu(F(z))\le1,
\quad F(0)=W_1,\quad F(1/2)=W_2.
\]
Set
\[
u=\frac1{\sqrt2}\begin{pmatrix}1\\i\end{pmatrix},
\quad f(z)=\langle F(z)u, u\rangle, \quad z\in\mathbb D.
\]
Then $f$ is a scalar Schur function with $f(0)=i/2,\quad f(1/2)=- i/2.$ By the Schwarz--Pick inequality,
\[
\frac45
=
\left|
\frac{f(0)-f(1/2)}
     {1-\overline{f(1/2)}f(0)}
\right|
\le\frac12,
\]
which is a contradiction. Therefore, the interpolation problem
is not solvable.
\end{example}
We end this section with a remark drawing connection to Sarason's approach to the Pick interpolation problem via the commutant lifting theorem.
\begin{remark} \label{rm:CLT}
    Recall that a bounded operator $T$ on a Hilbert space is a contraction if and only if $I - T^*T \succeq 0$. In classical operator-valued Nevanlinna--Pick interpolation, the positivity of $I - T^*T$ provides the canonical necessary and sufficient condition for solvability. More precisely, suppose an $N$-point interpolation problem $\lambda_j \mapsto W_j$ admits a holomorphic solution $F \colon \mathbb{D} \to M_d$ with $\|F\|_\infty = \sup_{z\in\mathbb{D}} \|F(z)\| \le 1$. The restriction of the contractive multiplication operator $M_F^*$ on the $\C^d$-valued Hardy-Hilbert space $H^2(\mathbb{C}^d)$ to the co-invariant subspace 
\[
\mathcal{M} = \operatorname{span}\{ k_{\lambda_j} \otimes v \colon j=1,\dots,N, \, v \in \mathbb{C}^d \} \subset H^2(\mathbb{C}^d),
\]
where for $\lambda \in \D$, $k_\lambda(z) \bydef (1-\bar{\lambda}z)$ is the Szeg\"o kernel function, defines a compressed operator $T = M_F^*|_{\mathcal{M}}$ acting via $T(k_{\lambda_j} \otimes v) = k_{\lambda_j} \otimes W_j^* v$. Consequently, $T$ is a contraction ($I - T^*T \succeq 0$), which is equivalent to the positive semi-definiteness of the classical operator Pick matrix
\[
\left[ \frac{I_d - W_i W_j^*}{1 - \lambda_i \bar{\lambda}_j} \right]_{i,j=1}^N \succeq 0.
\]
The sufficiency of this positivity condition is guaranteed by Sarason's Commutant Lifting Theorem: if $\mathcal{M} \subset H^2(\mathbb{C}^d)$ is a closed $M_z^*$-invariant subspace and $T \colon \mathcal{M} \to \mathcal{M}$ is a contraction commuting with $M_z^*|_{\mathcal{M}}$, then $T$ lifts to a contractive multiplication operator $M_F^*$ for some $F \in H^\infty(\D, M_d)$ with $\|F\|_\infty = \|T\| \le 1$, thereby ensuring interpolation by $F$.

In the proof of Theorem \ref{thm:Nec-cond}, we establishs that the positivity of the matrix family $P(\theta)$ for all $\theta \in [0, 2\pi)$ corresponding precisely to the compressed operator $T$ on $\mathcal{M}$ being a numerical contraction ($\nu(T) \le 1$). Were an analogue of Sarason's Commutant Lifting Theorem valid in this setting, any such $T$ would lift to a multiplication operator $M_F^*$ satisfying $\nu(M_F) = \nu(T) \le 1$ and $T=M_F^*|_{\mathcal M}$ for some $F\in H^\infty(\D, M_d)$. Applying \cite[Theorem~1]{Bebiano} to this lift would imply $\sup_{z \in \mathbb{D}} \nu(F(z)) = \nu(M_F) \leq 1$. Consequently, $F$ would be a solution of the $\nu$-Schur interpolation problem. However, Example \ref{exp:roadblock} reveals that the positivity conditions $P(\theta) \succeq 0$ are necessary but strictly insufficient, demonstrating that Sarason-type lifting fails for numerical contractions.
\end{remark}

\section{Structural factorizations of $\nu$-Schur functions} \label{sec:structure}
As Example~\ref{exp:roadblock} demonstrates, the natural necessary condition for the solvability of the interpolation problem fails to be sufficient, underscoring the difficulty of characterizing the solvability of a $\nu$-interpolation problem. To overcome this obstacle, we develop a function-theoretic tool that plays a pivotal role in establishing a necessary and sufficient solvability criterion in the next section.

We proceed with several preliminary results that are necessary for later developments.
\subsection{The $\nu$-Schur functions and their radial limits} \label{sub:basics}
Let $H^\infty(\D, M_d)$ be the Banach algebra of $M_d$-valued bounded holomorphic functions in $\D$ equipped with the sup-norm given by
$$
\|f\|_\infty \bydef \sup_{z\in \D} \|f(z)\|.
$$
We define another norm on $H^\infty(\D, M_d)$ given by
$$
 \|f\|_{\nu, \infty} \bydef \sup_{z\in\D} \nu(f(z)) \quad \text{for } f\in H^\infty(\D, M_d).
$$
It is evident from the well known inequality 
\begin{align}\label{eq:norm-equivalence}
\nu(T) \leq \|T\| \leq 2\nu(T),
\end{align}that for any $f\in H^\infty(\D, M_d)$,
\begin{align}\label{eq:fun-norm-equivalence}
    \|f\|_{\nu, \infty} \leq \|f\|_{\infty} \leq 2 \|f\|_{\nu, \infty}
\end{align}
proving the equivalence of the two norms above on $H^\infty(\D, M_d)$. It is important to note that unlike $\|\cdot\|_\infty$, the norm $\|\cdot\|_{\nu, \infty}$ does not make the Banach space $H^\infty(\D, M_d)$ a Banach algebra when $d>1$ because there are matrices $A,B$ so that $\nu(AB) > \nu(A)\nu(B)$; see for example \cite{Bhatia}.

Let us formally define the main holomorphic function theoretic thread of the paper.
\begin{definition}\label{def:NRC-holo-func}
A function $f$ in $H^\infty(\D, M_d)$ is said to be a {\it $\nu$-Schur function} if $\|f\|_{\nu, \infty} \leq 1$, a {\it strict $\nu$-Schur function} if $\|f\|_{\nu, \infty} < 1$. 
\end{definition}

Throughout this article we consider the standard normalized Lebesgue measure $m$ on $\T$. Equip $M_d$ with the Hilbert-Schmidt inner product
$$
\langle A, B\rangle_{\operatorname{HS}} \bydef \operatorname{tr}(B^*A) \quad 
\|A\|_{\operatorname{HS}} \bydef \operatorname{tr}(A^*A)^{1/2}.
$$With this inner product structure on $M_d$ we define
$$
L^2(\T, M_d) \bydef \left\lbrace \psi: \T \to M_d : \psi \text{ is measurable and } \|\psi\|^2_{\operatorname{HS},2} = \int_{\T} \|\psi(\zeta)\|_{\operatorname{HS}}^2 \ dm(\zeta) <\infty \right \rbrace.
$$ Note that $L^2(\T, M_d)\cong L^2\otimes M_d$ is a Hilbert space with the Hilbert-Schmidt inner product:
$$
\langle \psi_1, \psi_2 \rangle_{\operatorname{HS},2} = \int_\mathbb{T} \langle \psi_1(\zeta),\psi_2(\zeta) \rangle_{\operatorname{HS}}\;dm(\zeta) = \int_{\T} \operatorname{tr}(\psi_2(\zeta)^* \psi_1(\zeta)) \ dm(\zeta). 
$$With every notion of $L^2$ comes along a notion of $L^\infty$:
$$
L^\infty(\T, M_d) \bydef \left\lbrace \psi: \T \to M_d : \psi \text{ is measurable and } \|\psi\|_{\operatorname{ess},\infty} = \operatorname{ess \ sup}_{\zeta \in \T} \|\psi(\zeta)\| <\infty \right \rbrace.
$$

Since every norm on the finite dimensional space $M_d$ are equivalent, if we replace the Hilbert-Schmidt norm $\|\cdot \|_{\textup{HS}}$ by the numerical radius $\nu(\cdot)$ or the operator norm $\|\cdot \|$, we obtain the same Banach spaces with the following equivalent norms. 
For the purpose of this article sometimes we use the following equivalent $L^\infty$-norm:
$$
\|\psi\|_{\nu\operatorname{-ess}, \infty} \bydef \operatorname{ess \ sup}_{\zeta \in \T} \nu(\psi (\zeta)).  
$$

The second inequality in \eqref{eq:fun-norm-equivalence} implies that every $\nu$-Schur function $f$ is a function in $H^\infty(\D, M_d)$. By the Fatou's theorem $f$ has a non-tangential limit function $f^\flat$ in $L^\infty(\T, M_d)$ -- see e.g., \cite[Chapter~4]{RR-book}. In particular,
$$
f^\flat(\zeta) = \lim_{r \to 1^-} f(r\zeta) \ \text{ for almost every } \zeta \text{ in }\T,
$$
where the convergence is understood in operator norm and hence, in numerical radius norm on $M_d$. This readily shows that
$$
\|f^\flat\|_{\nu\operatorname{-ess}, \infty} \bydef \operatorname{ess \ sup}_{\zeta \in \T} \nu(f^\flat (\zeta)) \leq \|f\|_\infty.
$$
Moreover, we recover $f$ from $f^\flat$ via the Poisson integral formula:
\begin{align}\label{Poisson}
f(z) = \int_{\T} \frac{1-|z|^2}{|\zeta-z|^2} f^\flat(\zeta) \ dm(\zeta) \quad \text{for }z\in \D.
\end{align}
For any unit vector $x \in \C^d$, utilizing the fact that 
$
|\langle f^\flat(\zeta) x, x \rangle |\leq \|f^\flat\|_{\nu\operatorname{-ess}, \infty} $ for almost every $\zeta$ in $\T$,
the Poisson integral formula above along with continuity of the norm $\nu(\cdot)$ yields
$$
\nu(f(z)) \leq \|f^\flat\|_{\nu\operatorname{-ess}, \infty}.
$$
Hence, we have
\begin{align}\label{RadLim-Iso}
\|f\|_{\nu, \infty} = \|f^\flat\|_{\nu\operatorname{-ess}, \infty}.
\end{align}

For notational convenience, we sometimes use $f$ instead of $f^\flat$.

\subsection{Outer functions and the spectral density} 
 Let $H^2(\C^d)$ be the standard $\C^d$-valued Hardy-Hilbert space on $\D$. 
\begin{definition} 
    A function $F$ in $H^\infty(\D, M_d)$ is said to be {\it outer} if the range of the multiplication operator $M_F : H^2(\mathbb{C}^d) \rightarrow H^2(\mathbb{C}^d)$ defined by 
    $$
    M_F(\varphi)(z) \bydef F(z)\varphi(z) \quad \mbox{for } \varphi \in H^2(\mathbb{C}^d)
    $$
    is dense in $H^2(\mathbb{C}^d)$.
\end{definition}
Note that for an outer function $F$, its determinant $\det F(z)$ is non-vanishing in $\D$. Indeed, if $\det F(z_0)=0$ for some $z_0\in\D$, then there is $0\neq u\in \operatorname{Ker}F(z_0)^*$. Then for any function $g\in H^2(\C^d)$, 
we have 
$
0=\langle (Fg)(z_0), u\rangle.  
$
Since $F$ is outer, this is not possible unless $u=0$. In fact, what is known is that a function $F \in H^\infty(\D, M_d)$ is outer if and only if $\det F$ is a scalar-valued outer function. See the Lemma in page 125 of \cite{Helson-book}.

Positive semi-definite matrix-valued functions on $\mathbb T$ admitting a factorization of the form $\varphi^*\varphi$, where $\varphi$ is outer, play an important role in several areas of mathematics, including prediction theory, linear-quadratic control, and $H^\infty$ robust control. Such functions are commonly known as {\it spectral densities}. Part~\textup{(i)} of the following theorem characterizes those functions in $L^\infty(\mathbb T,M_d)$ that admit such a factorization and describes the corresponding outer spectral factors. This result originates in the work of Wiener and Masani \cite[Theorem~7.13]{Wiener-Masani-I} and Helson and Lowdenslager \cite[Theorem~9]{HL}; the formulation presented here follows Barclay \cite[Proposition~3.4]{Barclay-JLMS}. Part~\textup{(ii)} is its polynomial
refinement, given by the matrix Fej\'er--Riesz theorem
\cite[Section~6, page~416-417]{Rosenblatt}. 
These will serve as a key ingredient in establishing a factorization theorem for $\nu$-Schur functions.
\begin{theorem}\label{thm:spectral-factor}
Let $w\in L^\infty(\T,M_d)$ be positive semi-definite
almost everywhere.
\begin{enumerate}
\item[\textup{(i)}]
\textup{(Wiener--Masani; Helson--Lowdenslager
\cite{Wiener-Masani-I,HL}.)}
Then $\log\det w\in L^1(\T)$ if and only if there exist
outer functions $\Phi_R,\Phi_L\in H^\infty(\D,M_d)$
such that
\begin{align}\label{aux1}
 w(\zeta)
 =\Phi_R (\zeta)^*\Phi_R(\zeta)
 =\Phi_L(\zeta)\Phi_L(\zeta)^*
 \quad\text{for almost every }\zeta\in\T,
\end{align}and
\begin{align}\label{aux2}
\Phi_R(0) \succ 0 \quad \mbox{and} \quad \Phi_L(0) \succ 0.
\end{align}
These factors are uniquely determined by the condition \eqref{aux2}.

\item[\textup{(ii)}]
\textup{(Matrix Fej\'er--Riesz theorem
\cite{Rosenblatt}.)}
If, in addition,
\[
 w(\zeta)=\sum_{k=-m}^{m}w_k\zeta^k
 \quad\text{and}\quad
 w(\zeta)\succ0\quad(\zeta\in\T),
\]
then the normalized factors in \textup{(i)} are
matrix polynomials of degree at most $m$, and
\[
 \det\Phi_R(z)\ne0,\quad \det\Phi_L(z)\ne0
 \quad(z\in\overline{\D}).
\]
In this case, the factorization identities hold
at every point of $\T$.
\end{enumerate}
\end{theorem}
\begin{remark}\label{rem:bounded-inverse}
Suppose in addition that $w(\zeta)\ge\delta I_d$ a.e.\ for some $\delta>0$.
Then $$\Phi_R^{-1} \in H^\infty(\D,M_d) \quad \text{and} \quad \|\Phi_R^{-1}\|_\infty\le\delta^{-1/2}.$$
Indeed, for almost all \ $\zeta\in\T$ and every $x\in\mathbb{C}^d$,
$\|\Phi_R(\zeta)x\|^2=\langle w(\zeta)x,x\rangle\ge\delta\|x\|^2$. So, for any $\varphi\in H^2(\mathbb{C}^d)$,
\begin{align*}
   \|\Phi_R\varphi\|_2^2 = \int_{\T} \|\Phi_R(\zeta) \varphi(\zeta)\|^2\; dm(\zeta) \geq  \delta\|\varphi\|_2^2.
\end{align*}
Thus, $M_{\Phi_R}$ is bounded below and has closed range. Since $\Phi_R$ is outer, the range of $M_{\Phi_R}$ is also dense, hence $M_{\Phi_R}$ is invertible with $\|M_{\Phi_R}^{-1}\|\le\delta^{-1/2}$. The inverse commutes with
the shift $M_z$, so $M_{\Phi_R}^{-1}=M_\Psi$ for some
$\Psi\in H^\infty(\D,M_d)$. Then $M_{\Phi_R\Psi}=I_{H^2(\C^d)}$ which implies that $\Phi_R(z) \Psi(z)=I_d$. Hence,
$\Psi=\Phi_R^{-1}$. In a similar fashion, we can prove $\Phi_L^{-1} \in H^\infty(\D,M_d)$ and $\|\Phi_L^{-1}\|_\infty\le\delta^{-1/2}$.
\end{remark}

\subsection{Borel measurability of unique minimizers}
We record the following version of Borel selection theorem. It will be required to produce a Borel measurable choices of Ando completion guaranteed by Proposition \ref{prop:ando}.
\begin{lemma}\label{lem:Borel}
Let $\mathcal{X}$ and $\mathcal{Y}$ be compact metrizable spaces, let $\mathcal G\subset \mathcal{X}\times \mathcal{Y}$ be a closed set, and let $\mathcal{G}_x=\{y\in \mathcal{Y}:(x,y)\in \mathcal G\}$ for each $x\in \mathcal{X}$. Assume that $\mathcal G_x\neq\varnothing$ for every $x\in \mathcal{X}$. Let $\varphi:\mathcal{Y}\to\mathbb R$ be a continuous function such that for every $x\in \mathcal{X}$ the restriction $\varphi|_{\mathcal G_x}$ attains its minimum exactly at one point, denoted by $\Phi(x)$. Then $\Phi:\mathcal{X}\to \mathcal{Y}$ is Borel measurable.
\end{lemma}

\begin{proof}
Since $\mathcal G$ is closed, each fibre $\mathcal G_x$ is a closed subset of the compact space $\mathcal{Y}$, hence compact. So, $\varphi|_{\mathcal G_x}$ attains its minimum, and by the hypothesis it is attained exactly at one point.

For a closed set $C\subset \mathcal{Y}$ define an auxiliary function as
\[
m_C(x)\bydef \min\{\varphi(y):\ y\in G_x\cap C\}\in\mathbb R \cup\{+\infty\}, \quad x\in \mathcal{X},
\]
with the convention $\min\varnothing \bydef +\infty$.
The minimum is attained whenever $\mathcal G_x\cap C\ne\varnothing$, because $\mathcal G_x\cap C$ is compact. Write $m=m_\mathcal{Y}$. Then $m(x)=\varphi(\Phi(x)) \in \mathbb R$. Note that
\begin{equation}\label{eq:mCgem}
m_C(x)\ \ge\ m(x)\quad\text{for every } x\in \mathcal{X} \text{ and every closed } C\subset \mathcal{Y}.
\end{equation}

We now show that for any closed set $C \subset \mathcal{Y}$, $m_C$ is lower semi-continuous. Indeed, for $c\in \mathbb R$,
\begin{equation}\label{eq:sublevel}
\{x\in \mathcal{X}:\ m_C(x)\le c\}
=\Pi_\mathcal{X}\Big(\mathcal G\cap (\mathcal{X}\times C)\cap (\mathcal{X}\times\{y \in \mathcal{Y}:\varphi(y)\le c\})\Big),
\end{equation}
where $\Pi_\mathcal{X}$ denotes the projection from $\mathcal{X} \times \mathcal{Y}$ onto $\mathcal{X}$.
To verify the set-theoretic equality \eqref{eq:sublevel}, let $m_C(x)\le c<\infty$. Then $\mathcal G_x\cap C\neq\varnothing$ and the minimum is attained at some $y\in \mathcal G_x\cap C$ with $\varphi(y)\le c$, so $x$ lies in the right-hand side of \eqref{eq:sublevel}. Conversely, if $x$ lies in the right-hand side of \eqref{eq:sublevel}, then there is $y\in \mathcal G_x\cap C$ with $\varphi(y)\le c$, whence $m_C(x)\le c$.
 
The above equality implies that $\{x\in \mathcal{X}:\ m_C(x)\le c\}$ is the projection of the closed set $\mathcal G\cap (\mathcal{X}\times C)\cap (\mathcal{X}\times\{y \in \mathcal{Y}:\varphi(y)\le c\})$ and hence, compact. Thus, $m_C$ is lower semi-continuous, and in particular, Borel measurable as a map into $(-\infty,+\infty]$. 

Finally, we prove that $\Phi$ is Borel measurable by showing that for every closed set $C\subset \mathcal{Y}$, $\Phi^{-1}(C)$ is Borel.
For such a set $C$, utilizing the uniqueness of the minimizer $\Phi(x)$ we observe that 
\begin{align*}
  \Phi^{-1}(C) = \left \lbrace x\in \mathcal{X}: \ m_C(x) = m(x)\right \rbrace.
\end{align*}
Further, the function $m_C-m$ is Borel. Hence, $\Phi^{-1}(C)$ is the zero set of a Borel function and therefore Borel measurable.
\end{proof}

\subsection{Factorization of $\nu$-Schur functions into Schur-outer functions}
In this subsection, we prove a factorization result for $\nu$-Schur functions. Given a $\nu$-Schur function $f$, the key step is to apply Proposition~\ref{prop:ando} to the non-tangential limit values $f^\flat(\zeta)$ and obtain a suitable \textit{Ando completion} $Z(\zeta)$ so that the block matrix in \eqref{eq:ando-block} is positive semi-definite with $T=f^\flat(\zeta)$ and $Z=Z(\zeta)$. The preliminary results above come into play when making these suitable choices. Once the choices are made, the rest of the proof is a matter of appropriate assembling the components, which is where the following lemma comes handy.

\begin{lemma}\label{lem:completion-factorization}
Let $f\in H^\infty(\D,M_d)$, and let
$L\in L^\infty(\T,M_d)$ be self-adjoint. Suppose that
for some $\delta>0$,
\begin{equation}\label{Ando_completion}
 \begin{bmatrix}
  2I_d-L(\zeta) & f(\zeta)^*\\
  f(\zeta)   & L(\zeta)
 \end{bmatrix}
 \succeq \delta I_{2d}
 \ \text{ almost everywhere } \zeta\in\T.
\end{equation}
If $\Phi_1,\Phi_2\in H^\infty(\D,M_d)$ are outer such that
\[
 L=\Phi_1\Phi_1^*,\
 \text{ and }\
 2I_d-L=\Phi_2^*\Phi_2,
\]
then $G=\Phi_1^{-1}f\Phi_2^{-1}$ is a strict Schur function. In particular, $ f=\Phi_1G\Phi_2$ and
\[
 \Phi_1(\zeta)\Phi_1(\zeta)^*
       +\Phi_2(\zeta)^*\Phi_2(\zeta)=2I_d
\] for almost every $\zeta\in \T$.
\end{lemma}
\begin{proof}
Comparing the diagonal blocks in the inequality \eqref{Ando_completion}, we get 
\[
 \delta I_d\preceq L(\zeta)\preceq (2-\delta)I_d, \text{ and }\  \delta I_d\preceq 2I_d-L(\zeta)\preceq(2-\delta)I_d,
\]
for almost every $\zeta\in\mathbb T$. Since $\Phi_1$ and $\Phi_2$ are outer and $L(\zeta)\succeq \delta I_d$ almost everywhere, Remark~\ref{rem:bounded-inverse} therefore implies that their inverses $\Phi_1^{-1}$ and $\Phi_2^{-1}$ belong to
$H^\infty(\D,M_d)$. Thus, $G$ defines an element of $H^\infty(\D,M_d)$, and 
\[
f=\Phi_1G\Phi_2  \ \text{ and }\
 \Phi_1(\zeta)\Phi_1(\zeta)^*
       +\Phi_2(\zeta)^*\Phi_2(\zeta)=2I_d\
 \text{ for almost every } \zeta\in \T.
\]
It remains to show that $\|G\|_{\infty}<1$. For almost every $\zeta \in \mathbb{T}$, we write
\[
 C(\zeta)\bydef
 \begin{bmatrix}
  \Phi_2(\zeta)&0\\
  0&\Phi_1(\zeta)^*
 \end{bmatrix}.
\]
Then we note using the factorization hypothesis that
\begin{align*}
    \begin{bmatrix}
    I_d & G(\zeta)^*\\
    G(\zeta) & I_d
  \end{bmatrix} 
  &= C(\zeta)^{* -1}
  \begin{bmatrix}
    \Phi_2(\zeta)^{*} \Phi_2(\zeta) & f(\zeta)^*\\
    f(\zeta) &  \Phi_1(\zeta) \Phi_1(\zeta)^*
  \end{bmatrix}
  C(\zeta)^{-1}\\
  &=C(\zeta)^{* -1}
  \begin{bmatrix}
    2I_d - L(\zeta) & f(\zeta)^*\\
    f(\zeta) &  L(\zeta)
  \end{bmatrix}
  C(\zeta)^{-1}\\
  & \succeq \delta C(\zeta)^{* -1} C(\zeta)^{-1}.
\end{align*}
Moreover, 
\[
\Phi_1(\zeta)\Phi_1(\zeta)^*=L(\zeta)\preceq (2-\delta)I_d\ \text{ and } \Phi_2(\zeta)^*\Phi_2(\zeta)= 2I_d-L(\zeta)\preceq (2-\delta)I_d
\]
reveal that $\|C(\zeta)\|^2\le 2-\delta$ a.e., and hence, $\|C(\zeta)^{-1}\|^2\ge \frac {1}{(2-\delta)}$ a.e. Thus,
\[
 \begin{bmatrix}
  I_d & G(\zeta)^*\\
  G(\zeta) & I_d
 \end{bmatrix}\succeq \frac{\delta}{2-\delta}I_{2d} \quad \text{for almost every } \zeta \in \T. 
 \]
Consequently,
\[
 \|G(\zeta)\|
 \leq 1-\frac{\delta}{2-\delta}<1
 \quad \text{ for almost every } \zeta \in \T.
\]
This proves that $\|G\|_\infty<1$.
\end{proof}

The following lemma shows how the boundary condition propagates into the unit disk.
\begin{lemma}\label{L:GenForm}
Suppose $ \Phi_1, \Phi_2  \in H^\infty(\mathbb{D}, M_d)$ are such that 
\begin{align}\label{GenIneq}
\Phi_1(\zeta) \Phi_1(\zeta)^* + \Phi_2(\zeta)^* \Phi_2(\zeta) \preceq 2I_d \quad \text{for a.e. } \zeta \in \mathbb{T}.
\end{align} Then:
\begin{itemize}
    \item[(i)] The relation \eqref{GenIneq} continues to hold for every $z\in \D$, i.e.,
    \begin{align} \label{GenIneq-2}
        \Phi_1(z) \Phi_1(z)^* + \Phi_2(z)^* \Phi_2(z) \preceq 2I_d;
    \end{align}

    \item[(ii)] Given a $G\in H^\infty(\mathbb{D}, M_d)$, define the $M_d$-valued function
\begin{align}\label{GenForm}
f=\Phi_1G \Phi_2.
\end{align}Then $\|f\|_{\nu,\infty}\leq \|G\|_\infty$.
\end{itemize}
\end{lemma}
\begin{proof}
Let us recall a standard fact from the theory of Nevanlinna and Hardy classes of vector-valued functions \cite[Chapter~4]{RR-book}. For a space $(X, \|\cdot\|_X)$, any $X$-valued bounded holomorphic function $g: \D \to X$, satisfies 
\begin{align}\label{eq:submean-ineq}
    \|g(z)\|_X^p \leq \int_\T \frac{1-|z|^2}{|\zeta - z|^2}  \|g(\zeta)\|_X^p \ dm(\zeta).
\end{align}
for $0<p<\infty$. This follows from \cite[Inequality~4.19]{RR-book} (with the convex functions $\phi(t)= e^{pt}$).

To prove (i), let $x$ be a unit vector in $\C^d$. Consider function
$$
u_x(z) \bydef \|\Phi_1(z)^*x\|^2 + \|\Phi_2(z)x\|^2 \quad \text{for } z\in \D.
$$
Note that $\|\Phi_1(z)^*x\|^2= \|\Phi_1(z)^{\textup{T}} \bar{x}\|^2$. Applying \eqref{eq:submean-ineq} with $p=2$ to both the $\C^d$-valued bounded holomorphic functions: $z\mapsto \Phi_1(z)^{\textup{T}} \bar{x}$ and $z\mapsto \Phi_2(z)x$, we obtain
\begin{align*}
u_x(z) &\le \int_{\T} \frac{1-|z|^2}{|\zeta -z|^2} \left( \|\Phi_1(\zeta)^*x\|^2 + \|\Phi_2(\zeta)x\|^2\right) \ dm(\zeta)\\ 
&\le 2 \int_{\T} \frac{1-|z|^2}{|\zeta -z|^2} dm(\zeta) =2
\end{align*}
Since $x\in \C^d$ was arbitrary, the matrix inequality \eqref{GenIneq-2} follows.

For a proof of (ii), fix a unit vector $x\in\mathbb{C}^d$. Apply \eqref{eq:submean-ineq} with $p=1$ to the scalar function in $H^\infty(\D)$: $z\mapsto \langle f(z)x, x \rangle$, to get
\begin{align*}
    |\langle f(z)x, x \rangle| & \leq \int_{\T} \frac{1-|z|^2}{|\zeta -z|^2} |\langle f(\zeta) x, x \rangle| \ dm(\zeta) \\
    &=\int_{\T} \frac{1-|z|^2}{|\zeta -z|^2} |\langle G(\zeta)\Phi_2(\zeta)x, \Phi_1(\zeta)^*x \rangle| \ dm(\zeta) \\
    & \le \int_{\T} \frac{1-|z|^2}{|\zeta -z|^2} \frac{\|G\|_{\infty}}{2}\left(
       \|\Phi_2(\zeta)x\|^2+\| \Phi_1(\zeta)^*x \|^2
     \right) \ dm(\zeta) \le \|G\|_\infty.
\end{align*}
Since the unit vector $x$ was arbitrary, $\nu(f(z)) \leq \|G\|_\infty$ for every $z\in \D$.
\end{proof}
We show that every $\nu$-Schur function can be factorized as in \eqref{GenForm}. We establish this by first finding the factorization for strict $\nu$-Schur
functions, for which the factors have some additional features.
\begin{theorem} \label{thm:gen-structure}
Let $f \in H^\infty(\mathbb{D}, M_d)$ satisfy $\|f\|_{\nu,\infty} < 1$.
Then there exist outer functions $\Phi_1, \Phi_2 \in H^\infty(\mathbb{D}, M_d)$ and $G \in H^\infty(\mathbb{D}, M_d)$ with $\|G\|_\infty < 1$ satisfying:
\begin{align}
\label{required1} f(z) &= \Phi_1(z) G(z) \Phi_2(z) \quad \text{for } z \in \mathbb{D},\\
\label{required2} 
2I_d &= \Phi_1(\zeta) \Phi_1(\zeta)^* + \Phi_2(\zeta)^* \Phi_2(\zeta) \quad \text{for a.e. } \zeta \in \mathbb{T}.
\end{align}
\end{theorem}

\begin{proof}
Let $ 0 < c = \|f\|_{\nu,\infty} < 1$ and $F = \frac{1}{c} f$. Then $\|F\|_{\nu,\infty} = 1$. Thus, there exist a Borel measurable set $S \subset \mathbb{T}$ such that $m(\mathbb{T} \setminus S) = 0$ and for every $\zeta \in S$, $\lim_{r \to 1^-} F(r \zeta)$ exists and $\lim_{r \to 1^-} \nu(F(r\zeta)) \le 1$.
Let us choose $\widetilde{F}$, defined as follows, as a representative for the boundary value function of $F$.
\[
\widetilde{F}(\zeta) =
\begin{cases}
     \lim_{r \to 1^-} F(r \zeta)  & \text{ for } \zeta \in S; \\
     0 & \text{ for } \zeta \in \mathbb{T} \setminus S.
\end{cases}
\]
Therefore, for each $\zeta \in \mathbb{T}$, $\nu(\tilde{F}(\zeta)) \le 1$. Proposition~\ref{prop:ando} implies that for every $\zeta \in T$, the set
$$
\mathcal{S}(\zeta) \bydef \left \lbrace  Z \in M_d : \begin{bmatrix}
    2I_d - Z & \widetilde{F}(\zeta)^*\\
    \widetilde{F}(\zeta) & Z
  \end{bmatrix} \succeq 0  \right \rbrace
$$
is non-empty, bounded and convex subset of $M_d$. 
In fact, for every $\zeta \in \mathbb{T}$, 
$$\mathcal{S}(\zeta) \subseteq (M_d^{+})_2  \bydef \{Z \in M_d^{\rm sa} : \ 0 \preceq Z \preceq 2I_d\}.$$ 
Moreover, it is easy to see that $\mathcal{S}(\zeta)$ is closed in $M_d$ and hence, a compact subset of $(M_d^{+})_2$ for every $\zeta \in T$. To get a measurable choice of $Z\in \mathcal{S}(\zeta)$ by varying $\zeta$ over $\mathbb{T}$, we shall apply Lemma~\ref{lem:Borel}. Let us consider the continuous function
$$
\varphi:(M_d^{+})_2 \rightarrow \mathbb{R} \ \text{ such that } \ \varphi(T) \bydef \|T\|_{\textup{HS}}^2.
$$
The set 
$$
G= \left\lbrace (\zeta, Z) \in \mathbb{T} \times (M_d^{+})_2: \ Z \in \mathcal{S}(\zeta)  \right \rbrace
$$
is a closed subset of $\mathbb{T} \times (M_d^{+})_2$ as $\mathcal{S}(\zeta)$ is compact. 

We show that $\varphi|_{\mathcal{S}(\zeta)}$ attains its minimum exactly at one point. Indeed, the minimum value of the function $\varphi|_{\mathcal{S}(\zeta)}$ measures the squared distance of the zero matrix to the non-empty compact convex subset $\mathcal{S}(\zeta)$ of the Hilbert space $M_d$ endowed with the Hilbert-Schmidt inner product. Hence, by \cite[Chapter~I, Theorem~2.5]{Conway-book}, the squared distance is achieved by a unique element of $\mathcal{S}(\zeta)$.

Invoking Lemma~\ref{lem:Borel}, we infer that there exists a Borel measurable function $\widetilde{H}: \mathbb{T} \rightarrow (M_d^{+})_2$ such that 
\begin{align*}
    \begin{bmatrix}
    2I_d - \widetilde{H}(\zeta) & \widetilde{F}(\zeta)^*\\
    \widetilde{F}(\zeta) & \widetilde{H}(\zeta)
  \end{bmatrix} \succeq 0 \quad \text{ for every } \zeta \in \mathbb{T}. 
\end{align*}
This means 
\begin{align*}
    \begin{bmatrix}
    2I_d - \widetilde{H}(\zeta) & \frac{1}{c} f(\zeta)^*\\
    \frac{1}{c} f(\zeta) & \widetilde{H}(\zeta)
  \end{bmatrix} \succeq 0 \quad \text{ for every } \zeta \in \mathbb{T}. 
\end{align*}
Define, $H(\zeta) = c \widetilde{H}(\zeta) + (1-c)I_d$ for $\zeta \in \mathbb{T}$. Then for every $\zeta \in \mathbb{T}$,
\begin{align}\label{eq:lower-estimate}
    \begin{bmatrix}
    2I_d - H(\zeta) & f(\zeta)^*\\
    f(\zeta) & H(\zeta)
  \end{bmatrix} = c \begin{bmatrix}
    2I_d - \widetilde{H}(\zeta) & \frac{1}{c} f(\zeta)^*\\
    \frac{1}{c} f(\zeta) & \widetilde{H}(\zeta)
  \end{bmatrix} + (1-c)\begin{bmatrix}
    I_d & 0 \\ 0 & I_d  
  \end{bmatrix} \succeq (1-c)I_{2d}.
\end{align}
For every $\zeta \in \mathbb{T}$, the inequality $0 \preceq \widetilde{H}(\zeta) \preceq 2I_d$ yields: 
\begin{align*}
    (1-c)I_d \preceq H(\zeta) \preceq (1+c)I_d \quad \text{and} \\
    (1-c)I_d \preceq 2I_d - H(\zeta) \preceq (1+c)I_d.
\end{align*}
Thus, $H$ and $2I_d - H$ are both positive matrix-valued functions in $L^\infty(\mathbb{T}, M_d)$. Moreover,
$$
d \ \log(1-c) \le \  \int_{\mathbb{T}} \log(\det H(\zeta)) \ dm(\zeta) \ \le d \ \log(1+c)
$$
and the same upper and lower bounds hold for $\int_{\mathbb{T}} \log \det (2I_d - H(\zeta)) \ dm(\zeta)$. 
Now part~(i) of Theorem~\ref{thm:spectral-factor} provides us unique outer functions $\Phi_1$ and $\Phi_2$ in $H^\infty(\mathbb{D}, M_d)$ such that
\begin{align}\label{glue}
H(\zeta) = \Phi_1(\zeta) \Phi_1(\zeta)^* \text{ and } 2I_d - H(\zeta) = \Phi_2(\zeta)^* \Phi_2(\zeta) \ \text{ a.e. on } \mathbb{T},
\end{align}
and the matrices $\Phi_1(0)$ and $\Phi_2(0)$ are positive. Applying Lemma~\ref{lem:completion-factorization} with
$L=H$ and $\delta=1-c$, we get
$f=\Phi_1G\Phi_2$ for a Schur function $G$ with $\|G\|_\infty<1$, together with
the required boundary identity. This completes the proof.
\end{proof}

\begin{corollary} \label{cor:general-factorization}
If an $M_d$-valued function $f$ is a $\nu$-Schur function then it can be factorized as
$f = \Phi_1 G \Phi_2$,
where $\Phi_1, \Phi_2, G \in H^\infty(\mathbb{D}, M_d)$ are such that $\|G\|_\infty \le 1$ and
$\Phi_1(\zeta) \Phi_1(\zeta)^* + \Phi_2(\zeta)^* \Phi_2(\zeta) \preceq 2I_d$ for almost every $\zeta \in \mathbb{T}$.
\end{corollary}
\begin{proof}
Let $0<r_n<1$ be a sequence increasing to $1$. Consider the function $f_n(z)=r_nf(r_nz)$. Then each $f_n$ is holomorphic in a neighbourhood of $\overline{\mathbb{D}}$ and $\|f_n\|_{\nu,\infty}\leq r_n<1$. By Theorem \ref{thm:gen-structure}, there are sequences $\{\Phi_1^{(n)}\}_n$ and $\{\Phi_2^{(n)}\}_n$ of outer functions and $\{G_n\}_n$ of strict Schur functions so that $f_n=\Phi_1^{(n)}G_n\Phi_2^{(n)}$ and
$$
\Phi_1^{(n)}(\zeta)\Phi_1^{(n)}(\zeta)^* + \Phi_2^{(n)}(\zeta)^* \Phi_2^{(n)}(\zeta) = 2 I_d \quad \text{for a.e. } \zeta\in\mathbb{T}.
$$
Then from Lemma~\ref{L:GenForm}, we get
\begin{align}\label{eq:ineq-inside}
\Phi_1^{(n)}(z)\Phi_1^{(n)}(z)^* + \Phi_2^{(n)}(z)^* \Phi_2^{(n)}(z) \preceq 2 I_d \quad \text{for } z \in\D.
\end{align}
We now employ Montel's theorem, to obtain three subsequences $\{G_{n_k}\}_k$, $\{\Phi_1^{(n_k)}\}_k$ and $\{\Phi_1^{(n_k)}\}_k$ converging in the compact-open topology to a Schur function $G$, and $H^\infty$ functions $\Phi_1$ and $\Phi_2$, respectively. It is clear that we have the factorization $f=\Phi_1 G \Phi_2$. 
Also, taking pointwise limit in the inequality \eqref{eq:ineq-inside} yields same inequality involving $\Phi_1$ and $\Phi_2$. Thus, considering the radial limits of $\Phi_1$ and $\Phi_2$ we see that
they satisfy the boundary condition almost everywhere on $\T$ as stated.
\end{proof}

We are able to achieve outer functions $\Phi_1,\Phi_2$ in Theorem \ref{thm:gen-structure} but not in Corollary \ref{cor:general-factorization}. The functions $\Phi_1$ and $\Phi_2$ in Corollary~3.9 may not, in general, be outer. We demonstrate this by a concrete example.
\begin{example}
\label{ex:nonouter-factorization}
Consider the constant function
\[
    f(z)=
    \begin{bmatrix}
        0&2\\
        0&0
    \end{bmatrix},
    \quad z\in\D.
\]
Since the numerical radius of the above matrix is one, $\|f\|_{\nu,\infty}=1$.
Suppose that
\[
    f=\Phi_1G\Phi_2,
    \quad
    \Phi_1,\Phi_2,G\in H^\infty(\D,M_2),
    \quad
    \|G\|_\infty\le1,
\]
and
\begin{equation}
\label{eq:nonouter-energy}
    \Phi_1(\zeta)\Phi_1(\zeta)^*
    +\Phi_2(\zeta)^*\Phi_2(\zeta)
    \preceq 2I_2
    \quad\text{for a.e. }\zeta\in\T.
\end{equation}
Here, boundary values are understood as radial limits.
We show that neither $\Phi_1$ nor $\Phi_2$ can be outer.
Fix a boundary point $\zeta\in\mathbb T$ at which the radial limits and the
preceding relations hold, and write
\[
    A=\Phi_1(\zeta)\Phi_1(\zeta)^*,
    \quad
    B=\Phi_2(\zeta)^*\Phi_2(\zeta).
\]
Then $A,B\succeq0$ and $A+B\preceq2I_2$. In particular,
$A_{11}\le2$ and $B_{22}\le2$. If $e_1,e_2$ denote the
standard basis vectors of $\C^2$, then
\begin{align*}
    2=|\langle f(\zeta)e_2,e_1\rangle|
    \le
    \|\Phi_2(\zeta)e_2\|\,
    \|\Phi_1(\zeta)^*e_1\|
    =\sqrt{B_{22}A_{11}}
    \le2.
\end{align*}
Consequently, $A_{11}=B_{22}=2$. Using $A+B\preceq2I_2$, we get that $A_{22}=B_{11}=0$. 

Since $A$ and $B$ are positive semidefinite, their
off-diagonal entries also vanish. Therefore,
\[
    \Phi_1(\zeta)\Phi_1(\zeta)^*
    =
    \begin{bmatrix}
        2&0\\
        0&0
    \end{bmatrix},
    \quad
    \Phi_2(\zeta)^*\Phi_2(\zeta)
    =
    \begin{bmatrix}
        0&0\\
        0&2
    \end{bmatrix}
    \quad\text{for a.e. }\zeta\in\T.
\]

In particular,
\[
    \|\Phi_1(\zeta)^*e_2\|^2=0,
    \quad
    \|\Phi_2(\zeta)e_1\|^2=0
    \quad\text{for a.e. }\zeta\in\T.
\]
Thus the second row of $\Phi_1$ and the first column of
$\Phi_2$ have zero boundary values almost everywhere.
Since each of the co-ordinate entries of $\Phi_1$ and $\Phi_2$ belongs to $H^\infty(\D)$, the second row of $\Phi_1$ and first column of $\Phi_2$ vanish identically. It follows that
\[
    \det\Phi_1(z)=\det\Phi_2(z)=0
    \quad \text{for }z\in\D.
\]
 Hence, neither $\Phi_1$ nor
$\Phi_2$ is outer.

Nevertheless, the factorization required by the corollary
does exist. For example, take the constant functions
\[
    \Phi_1=
    \begin{bmatrix}
        \sqrt{2}&0\\
        0&0
    \end{bmatrix},
    \quad
    G=
    \begin{bmatrix}
        0&1\\
        0&0
    \end{bmatrix},
    \ \text{ and }
    \Phi_2=
    \begin{bmatrix}
        0&0\\
        0&\sqrt{2}
    \end{bmatrix}.
\]
Then
\[
    f=\Phi_1G\Phi_2,
    \quad
    \|G\|_\infty=1,
    \ \text{ and }
    \Phi_1\Phi_1^*+\Phi_2^*\Phi_2=2I_2.
\]
Thus, even when the boundary condition holds with
equality, the two square factors may fail
to be outer.
\end{example}

We next show that when the strict $\nu$-Schur function in Theorem \ref{thm:gen-structure} is continuous up to the boundary $\mathbb{T}$, then the outer functions $\Phi_1,\Phi_2$ in it can be taken to be {\em outer polynomials}.

\begin{theorem} \label{thm:strict-factorization}
Let $f:\mathbb{D}\to M_d(\mathbb{C})$ be a strict $\nu$-Schur function (i.e., $\|f\|_{\nu,\infty} < 1$) that is continuous up to the boundary $\mathbb{T}$. Then there exist two outer polynomials $P,Q$ and a strict Schur function $G$ (i.e., $\|G\|_{\infty}<1$)
such that $f$ gets factored as
\begin{align}\label{TheFactor}
   f = PGQ
\end{align}  
and $P,Q$ satisfy the boundary condition 
\begin{align}\label{BoundaryCond}
P(\zeta)P(\zeta)^* + Q(\zeta)^*Q(\zeta) = 2I_d \quad\mbox{for all }\zeta \in \mathbb{T}.
\end{align} Furthermore, $G$ can be chosen to be rational if $f$ is so. Moreover, $P$ and $Q$ can be chosen to be invertible on $\overline{\mathbb{D}}$.
\end{theorem}
\begin{proof}
The advantage of $f$ being a strict $\nu$-Schur function is that there is a $0<c<1$ such that
\begin{equation}\label{eq:r-choice}
  \|f\|_{\nu,\infty} = \sup_{\zeta\in\mathbb{T}}\nu(f(\zeta))<c<1.
\end{equation}
We shall take this advantage and apply Proposition~\ref{prop:ando} to
$f(\zeta)/c$ instead of $f(\zeta)$, for each $\zeta\in\mathbb{T}$, to obtain a self-adjoint matrix $Z_\zeta$ such that
\begin{equation}\label{eq:pointwise-ando}
  \begin{bmatrix}
    2I_d - Z_\zeta & f(\zeta)^*/c\\
    f(\zeta)/c & Z_\zeta
  \end{bmatrix}\succeq 0.
\end{equation}
Setting $\widetilde Z_\zeta \bydef cZ_\zeta+(1-c)I_d$, we observe from \eqref{eq:pointwise-ando} that
\begin{align} 
 \begin{bmatrix}
   2I_d - \widetilde Z_\zeta & f(\zeta)^*\\
   f(\zeta)& \widetilde Z_\zeta
 \end{bmatrix}
 &= c
 \begin{bmatrix}
   2I_d - Z_\zeta & f(\zeta)^*/c\\
   f(\zeta)/c & Z_\zeta
 \end{bmatrix}
 +(1-c)I_{2d}\notag\\
 &\succeq(1-c)I_{2d}.
 \label{eq:strict-pointwise-block}
\end{align}

Now, $f$ being continuous on $\mathbb{T}$, for every fixed $\zeta \in\mathbb{T}$, the map
$$
\mathbb{T}\ni \eta \mapsto \begin{bmatrix}
   2I_d - \widetilde Z_\zeta & f(\eta)^*\\
   f(\eta)& \widetilde Z_\zeta 
 \end{bmatrix}\in M_d^{\rm sa}
$$
is continuous. In particular, continuity of this map at $\zeta\in\mathbb{T}$ along with \eqref{eq:strict-pointwise-block} implies that there is an open neighborhood $U_\zeta$ of $\zeta$ in $\mathbb{T}$ such that
\begin{equation}\label{eq:local-fixed-A}
  \begin{bmatrix}
   2I_d - \widetilde Z_\zeta & f(\eta)^*\\
   f(\eta) & \widetilde Z_\zeta
  \end{bmatrix}
  \succ \frac{1-c}{2}I_{2d},
  \quad \text{for every } \eta\in U_\zeta.
\end{equation}
Compactness of $\mathbb{T}$ yields a finite subcover, say, $U_{\zeta_1},\ldots,U_{\zeta_s}$ of $\mathbb{T}$. We choose a continuous partition of unity
$\psi_1,\ldots,\psi_s$ subordinate to this subcover, i.e., 
\begin{align*}
0\leq \psi_j\leq 1 \mbox{ for each } j,\quad \sum_j\psi_j=1, \quad \mbox{and}\quad 
\operatorname{supp}\psi_j\subseteq U_{\zeta_j}.
\end{align*}

Consider the function $\Psi$ on $\mathbb{T}$ defined as
$$
  \Psi(\zeta)\bydef \sum_{j=1}^s\psi_j(\zeta)\widetilde Z_{\zeta_j}, \ \, \zeta\in \mathbb{T}.
$$
Clearly, $\Psi$ is a continuous self-adjoint matrix-valued function. So, using uniform approximation of continuous functions on $\mathbb{T}$ by Laurent polynomials, applied entry-wise and followed by a self-adjoint symmetrization, we can find a self-adjoint matrix-valued Laurent polynomial $L$ such that
\begin{align} \label{eq:trig-approx}
   \sup_{\zeta\in\mathbb{T}}\|L(\zeta)- \Psi(\zeta)\| <\frac{1-c}{4}. 
\end{align}

Further, 
the properties of the partition of unity $\{\psi_j\}_j$ and equation \eqref{eq:local-fixed-A} imply that for each $\zeta\in\mathbb{T}$,
\begin{align}
  \begin{bmatrix}
    2I_d- \Psi(\zeta) & f(\zeta)^*\\
    f(\zeta) & \Psi(\zeta)
  \end{bmatrix} & =
  \begin{bmatrix}
      \sum_{j=1}^s\psi_j(\zeta) (2I_d - \widetilde{Z}_{\zeta_j}) & \sum_{j=1}^s\psi_j(\zeta) f(\zeta)^*\\
    \sum_{j=1}^s\psi_j(\zeta) f(\zeta) & \sum_{j=1}^s\psi_j(\zeta)\widetilde Z_{\zeta_j}
  \end{bmatrix}  \notag \\
  &= \sum_{j=1}^s \psi_j(\zeta) \begin{bmatrix}
2I_d - \widetilde{Z}_{\zeta_j} &f(\zeta)^*\\
    f(\zeta) & \widetilde Z_{\zeta_j}
  \end{bmatrix} 
  \succ \frac{1-c}{2}I_{2d}. \label{eq:continuous-strict-block}
\end{align}  
Further more, for each $\zeta \in \mathbb{T}$, we also have
\begin{align}
  \begin{bmatrix}
    2I_d- L(\zeta) & f(\zeta)^*\\
    f(\zeta) & L(\zeta)
  \end{bmatrix} &= \begin{bmatrix}
    2I_d - \Psi(\zeta) & f(\zeta)^*\\
    f(\zeta) & \Psi(\zeta)
  \end{bmatrix} + \begin{bmatrix}
    \Psi(\zeta) - L(\zeta) &0 \\0& L(\zeta)-\Psi(\zeta)
  \end{bmatrix}  \notag\\
  & \succeq \frac{1-c}{4}I_{2d}. \label{eq:Laurent-strict-block}
\end{align}
The penultimate inequality is obtained from \eqref{eq:trig-approx} and \eqref{eq:continuous-strict-block}.
Thus, from \eqref{eq:Laurent-strict-block} we see that on $\mathbb{T}$
$$
  L \succeq \frac{1-c}{4}I_d,\quad \text{and} \quad 
  2I_d-L \succeq \frac{1-c}{4}I_d.
$$
Applying part~(ii) of Theorem~\ref{thm:spectral-factor} in the left-factor form to $L$ and in the
right-factor form to $2I_d - L$, we obtain outer polynomials $P,Q$ such that
\begin{equation}\label{eq:spectral-factors}
  L(\zeta)=P(\zeta)P(\zeta)^*\quad \mbox{and}\quad
  2I_d - L(\zeta) = Q(\zeta)^* Q(\zeta)
\end{equation}for each $\zeta\in\mathbb{T}$. 
Applying Lemma~\ref{lem:completion-factorization},  with
$\Phi_1=P$, $\Phi_2=Q$, and $\delta=(1-c)/4$, we get $f=PGQ$ with $\|G\|_\infty<1$ and the required
boundary identity. 

Since
$G=P^{-1}fQ^{-1}$, the function $G$ is rational whenever
$f$ is rational. Moreover, the strict positivity of
$L$ and $2I_d-L$ on $\T$, together with the invertibility
of the outer factors in $\D$, shows that $P$ and $Q$
are invertible on $\overline{\D}$.
\end{proof}

Let the functions $f, G, \Phi_1$ and $\Phi_2$ be as in Corollary \ref{cor:general-factorization}. The boundary condition on $\Phi_j$ as in Corollary \ref{cor:general-factorization} ensures that with $\Psi_j=\frac{1}{\sqrt{2}}\Phi_j$ for each $j=1,2$, we have
\begin{align}\label{boundaryInEq}
\Psi_1(\zeta)\Psi_1(\zeta)^* + \Psi_2(\zeta)^*\Psi_2(\zeta) \preceq I_d
\end{align}almost everywhere on $\mathbb{T}$. Consider the rectangular matrix-valued function defined on $\mathbb{D}$ as follows:
\begin{align}\label{V}
V(z)=\begin{bmatrix}
        \Psi_1(z)^* \\
        \Psi_2(z)
    \end{bmatrix}: \mathbb{C}^d \mapsto \begin{bmatrix}
        \mathbb{C}^d \\
        \mathbb{C}^d
    \end{bmatrix}.
\end{align}
Then using \eqref{boundaryInEq} and part (i) of Lemma~\ref{GenForm}, we obtain
\begin{align*}
V(z)^*V(z) = \Psi_1(z)\Psi_1(z)^* + \Psi_2(z)^*\Psi_2(z) \preceq I_d \quad \text{for } z\in \D, 
\end{align*}
and thereby, $\|V\|_\infty \leq 1$.

Note that equality in \eqref{boundaryInEq} would imply that $V(\zeta)$ is isometry valued almost everywhere on $\mathbb{T}$.

Further, with $\Psi_1,\Psi_2$ and $V$ as above, the identity $f=\Phi_1 G \Phi_2$ as in Corollary \ref{cor:general-factorization} can be expressed as
\begin{align}\label{DilationYes!}
        f =  V^*\begin{bmatrix}
         0 & 2G\\
         0 & 0
    \end{bmatrix}V.
\end{align}
With this observations in place, we now summarize the structural results obtained in this section in the theorem below:
\begin{theorem}\label{T:StructureSummary}
    Let $f$ be an $M_d$-valued function. Then the following are equivalent:
    \begin{enumerate}
        \item[\bf({A)}] $f$ is a $\nu$-Schur function;
        \item[\bf({B)}] There exist Schur functions $G, \Psi_1$ and $\Psi_2$ so that $V$ as in \eqref{V} is contractive function on $\mathbb{D}$ and $f$ can be expressed as \eqref{DilationYes!}.
    \end{enumerate}
    Moreover, the functions $G,\Psi_1,\Psi_2$ can be chosen with the following bonuses:
    \begin{enumerate}
        \item[\bf({C)}] if $\|f\|_{\nu,\infty}<1$, then $V$ can be chosen to take isometry valued almost everywhere on the unit circle $\mathbb{T}$;
        \item[\bf({D)}] if $f$ is continuous up to the boundary $\mathbb{T}$ and $\|f\|_{\nu,\infty}<1$, then in addition to {\bf({C)}}, $\Psi_1,\Psi_2$ can be chosen to be polynomials; and
        \item[\bf({E)}] if $f$ is rational with no poles in $\overline{\mathbb{D}}$ and $\|f\|_{\nu,\infty}<1$, then in addition to {\bf({D)}}, $G$ can be chosen to be rational as well.
    \end{enumerate}
\end{theorem}
\begin{proof}
    Equivalence of {\bf(A)} and {\bf(B)} is given by Lemma \ref{L:GenForm} and Corollary \ref{cor:general-factorization}; {\bf(C)} is the assertion of Theorem \ref{thm:gen-structure} while {\bf(D)} is proved in Theorem \ref{thm:strict-factorization}.
\end{proof}

\section{Solvability criteria for numerical radius interpolation problem} \label{sec:interpolation}

In this section, we present two distinct solvability criteria for the numerical radius interpolation problem, each detailed in its own dedicated subsection.
\subsection{Solvability criterion via structural factorizations of $\nu$-Schur functions}
We first present a solvability criterion using the factorization result of the $\nu$-Schur functions. 

The solvability criterion for the classical matrix-valued Nevanlinna--Pick interpolation problem will be used -- see \cite{AM-book}. 

\begin{theorem}\label{thm:matrix-NP}
Let $\lambda_1, \ldots, \lambda_N \in \mathbb{D}$ be distinct nodes and $G_1, \ldots, G_N \in M_d$. There exists $G \in H^\infty(\mathbb{D}, M_d)$ satisfying $\|G\|_\infty \le 1$ and $G(\lambda_j) = G_j$ ($1 \le j \le N$) if and only if 
\begin{align*}
 \left[
      \frac{I_d - G_j G_k^*}{1 - \lambda_j \bar{\lambda}_k}
    \right]_{j,k=1}^N \succeq 0.
\end{align*}
\end{theorem}
A solvable Pick interpolation problem is said to be extremal, if there is no solution $F$ of strict Schur norm $\|F\|_\infty <1$; see Chapter~6 of \cite{AM-book}. We now introduce a parallel terminology in the present setting.
\begin{definition}
    A solvable $\nu$-interpolation problem is said to be {\it extremal}, if every solution is of  $\nu$-Schur norm $\|\cdot\|_{\nu,\infty}$ one.
\end{definition}
A solvable $\nu$-interpolation problem is therefore non-extremal, if it admits a solution $F$ satisfying $\|F\|_{\nu,\infty} < 1$. Consequently, every solvable $\nu$-interpolation problem is either extremal or non-extremal. We now state and prove our interpolation theorem.
\begin{theorem} \label{thm:interpolation}
The interpolation problem
\begin{align}\label{Solvability}
\bm\Lambda=(\lambda_1, \ldots, \lambda_N) \mapsto (W_1, \ldots, W_N)=\textbf{W}
\end{align} is solvable by a $\nu$-Schur function if and only if
there exist two functions $\Phi_1, \Phi_2$ in $H^\infty(\mathbb D, M_d)$ and $d\times d$ matrices $G_1, \ldots, G_N $ such that
\begin{enumerate}
\item[\bf{(NS1)}] $W_j = \Phi_1(\lambda_j) G_j \Phi_2(\lambda_j)$ for every $j=1,2,\dots,N$;
\item[\bf{(NS2)}] the Pick matrix
  \begin{equation}\label{eq:classical-pick}
    \left[
      \frac{I_d - G_j G_k^*}{1 - \lambda_j \bar{\lambda}_k}
    \right]_{j,k=1}^N \succeq 0; \quad \mbox{and}
  \end{equation}
\item[\bf{(NS3)}] the boundary condition holds:
  \begin{equation}\label{eq:energy-inequality}
    \Phi_1(\zeta) \Phi_1(\zeta)^* + \Phi_2(\zeta)^* \Phi_2(\zeta) \preceq 2I_d \quad \text{for } \zeta \in \mathbb{T}.
  \end{equation}
\end{enumerate}
Furthermore, we have the following bonuses: 
\begin{enumerate}
    \item[\bf{(NE)}] If the problem \eqref{Solvability} is non-extremal, i.e., if it is solvable by a function $F$ with $\|F\|_{\nu,\infty}<1$, then the matrix in {\bf{(NS2)}} is positive definite and the functions $\Phi_1$ and $\Phi_2$ can be chosen to be outer polynomial and to have equality in {\bf{(NS3)}};
    \item[\bf{(E)}] If the problem is extremal, then there is a sequence of target data $\textbf{W}_n$ so that the problem $\bm\Lambda\mapsto \textbf{W}_n$ is non-extremal and $\textbf{W}_n\to \textbf{W}$.
\end{enumerate}
\end{theorem}
\begin{proof}
    Suppose the problem \eqref{Solvability} is solvable by the $\nu$-Schur function $F$. Apply Corollary \ref{cor:general-factorization} to $F$ to obtain a Schur function $G$ and outer functions $\Phi_1,\Phi_2$ in $H^\infty(\mathbb{D}, M_d)$ so that $F=\Phi_1G\Phi_2$ and for almost every $\zeta$ in $\mathbb{T}$,
    $$
    \Phi_1(\zeta)\Phi_1(\zeta)^* + \Phi_2(\zeta)^*\Phi_2(\zeta) \preceq 2I_d.
    $$Choose $G_j=G(\lambda_j)$ for each $j=1,2,\dots, N$. Condition (\textbf{NS1}) and (\textbf{NS3}) are therefore readily achieved. By the classical Nevanlinna Pick Interpolation Theorem \ref{thm:matrix-NP}, condition (\textbf{NS2}) is also satisfied.

    Conversely, suppose there are outer functions $\Phi_1, \Phi_2 \in H^\infty(\mathbb D, M_d)$ and $d\times d$ matrices $G_1, \ldots, G_N $ such that conditions (\textbf{NS1})-(\textbf{NS3}) hold. By the classical Nevanlinna Pick Interpolation Theorem \ref{thm:matrix-NP} applied to condition (\textbf{NS2}), there is a Schur function $G$ that interpolates the data. By condition (\textbf{NS1}), the function
    $$
    F=\Phi_1G\Phi_2
    $$also interpolates the data. By Lemma \ref{L:GenForm}, $\|F\|_{\nu,\infty}\leq \|G\|_\infty \leq 1.$
\smallskip

\noindent 
{\sf{Proof of {\bf{(NE)}}}:} Let us suppose that the interpolation problem is non-extremal, i.e., there is $0<c<1$ and 
a solution $F\in H^\infty(\mathbb{D},M_d)$ satisfying $\|F\|_{\nu,\infty} < c$. Let $0<s<1$ and consider the dilation $F_s(z)=F(sz)$ of $F$. The function $F_s$ is analytic on the
disk of radius $1/s>1$, and
\begin{equation}\label{eq:Fs-bound}
  \sup_{z\in\overline{\mathbb{D}}}\nu(F_s(z)) \le \|F\|_{\nu,\infty} < c.
\end{equation}
Moreover, for each fixed $j$,
\[
  F_s(\lambda_j)=F(s\lambda_j)\longrightarrow F(\lambda_j)=W_j
  \quad \text{as }s \rightarrow 1^-.
\]

Let $\ell_1,\ldots,\ell_N$ be the Lagrange polynomials given by
\begin{align}\label{LagrangePoly}
    \ell_j(z) = \prod_{k \neq j} \frac{z -\lambda_k}{\lambda_j -\lambda_k} \ \text{ for } j=1,\ldots, N.
\end{align} These polynomials have the property:
\begin{align}\label{LagId}
\ell_j(\lambda_i) = \delta_{i,j} = \begin{cases}
    1 \quad \mbox{if } i=j\\
    0 \quad \mbox{otherwise}.
\end{cases}
\end{align}This is to be used in what follows. It is convenient to denote the tuple $(\lambda_1,\dots,\lambda_N)$ by $\bm\Lambda$ and set
\[
  C_{\bf\Lambda}\bydef\sum_{i=1}^N\sup_{z\in\overline{\mathbb{D}}}|\ell_i(z)|<\infty.
\]
We choose $s$ sufficiently close to $1$ so that
\begin{equation}\label{eq:node-error-small}
  C_{\bf\Lambda} \max_j \|W_j-F_s(\lambda_j)\|
  <\frac{1-c}{2}.
\end{equation}
Define the $M_d$-valued function
\begin{equation}\label{eq:exact-polynomial-interpolant}
  H(z)=F_s(z)+\sum_{j=1}^N \ell_j(z) \bigl(W_j - F(\lambda_j)\bigr).
\end{equation}
Clearly, $H$ is holomorphic in a neighbourhood of $\overline{\D}$. The Lagrange identities \eqref{LagId} give $H(\lambda_j)=W_j$ for every $j$. Furthermore, for $z\in\overline{\mathbb{D}}$,
\begin{align*}
 \|H(z)-F_s(z)\|
 &\le
   \sum_{j=1}^N|\ell_j(z)|\,
       \|W_j-F_s(\lambda_j)\|\\
 &\le C_{\bf\Lambda}\max_j \|W_j -F_s(\lambda_j)\| <\frac{1-c}{2},
\end{align*}
where to obtain the last inequality we used \eqref{eq:node-error-small}. Combining this
with \eqref{eq:Fs-bound} and the triangle inequality for the numerical radius
gives
\[
  \|H\|_{\nu,\infty} 
  \le \|F_s\|_{\nu,\infty} + \|H - F_s\|_{\nu,\infty}
  < c + 2 \|H - F_s\|_{\infty}  < c + (1-c) =1.
\]
Theorem~\ref{thm:strict-factorization} applied to $H$ gives
\[
  H=PGQ,\quad \|G\|_\infty <1,\quad
  PP^*+Q^*Q=2I_d\ \text{ on }\mathbb{T},
\]
for two outer polynomials $P,Q$ in $M_d(\mathbb{C}[z])$ and a strict Schur function $G$. Choose $G_j=G(\lambda_j)$ for $j=1, \ldots, N$, and $(\Phi_1,\Phi_2)=(P,Q)$. Since $G$
is a strictly contractive function, Theorem~\ref{thm:matrix-NP} shows that the Pick matrix \eqref{eq:classical-pick} in {\bf{(NS2)}} is positive definite. As $H(\lambda_j)=W_j$, for each $j=1,2,\dots, N$, we have completed the proof of {\bf{(NE)}}.



\smallskip

\noindent 
{\sf{Proof of {\bf{(E)}}}:} Suppose the problem \eqref{Solvability} is extremal and $F$ is a solution with $\|F\|_{\nu, \infty} =1$. Given a sequence $\{0<r_n<1\}_n$ increasing to $1$ as $n\to \infty$, define
$$
F_n(z) \bydef r_n F(r_n z) \ \text{ for } z\in \mathbb{D}.
$$
Each $F_n$ is holomorphic in the open disk of radius $1/r_n$ around $0$ and $\|F_n \|_{\nu, \infty} \le r_n $. 

For each $n\geq 1$, let us set
$$
\boldsymbol{W}^{(n)} = (F_n(\lambda_1), \dots, F_n(\lambda_N)).
$$

It is easy to observe that $F_n$ converges to $F$ in the compact-open topology as $n$ tends to $\infty$. In particular, the convergence 
$\boldsymbol{W}^{(n)} \to \boldsymbol{W}$ follows.
\end{proof}

\subsection{Solvability criterion via a Nehari-type distance formula}\label{sub:Nehari-type formula}

The classical Nehari theory identifies the distance of a scalar $L^\infty(\T)$ symbol $\varphi$ from $H^\infty(\D)$ with the norm of the Hankel operator $H_\varphi$: $$ \operatorname{dist}_{L^\infty(\T)}(\varphi, H^\infty(\D)) = \|(I-P_{H^2(\D)})M_\varphi|_{H^2(\D)}\|,$$
where $M_\varphi$ denotes the multiplication by $\varphi$ on $L^2(\T)$ and $P_{H^2(\D)}$ is the orthogonal projection from $L^2(\T)$ onto $H^2(\D)$.
The classical holomorphic interpolation problems are intricately connected with the Nehari theory as well as model spaces \cite{FFGK, Peller-book}. We exploit these ideas to provide a distance criterion for our numerical radius interpolation problem. 

For the purpose of this subsection, we define a new norm on the finite dimensional spaces $M_d$ via its dual space as follows:
$$
\nu^*(T) \bydef \sup \{ |\operatorname{tr}(A^*T)|: A\in M_d \text{ with } \nu(A) \leq 1\}.
$$

We define $M_d$-valued $L^1$ space as
$$
L^1_{\nu^*}(\T,M_d) \bydef \left\lbrace f:\T \to M_d: \ f \text{ is measurable and } \|f\|_{\nu^*, 1} \bydef \int_{\T} \nu^*(f(\zeta)) \ dm(\zeta) <\infty \right\rbrace,
$$
equipped with the $\|\cdot \|_{\nu^*, 1}$ norm. We equip the space $L^\infty(\T, M_d)$ with the norm 
$\|\cdot\|_{\nu\operatorname{-ess}, \infty}$ as in Subsection~\ref{sub:basics} and denote it by $L^\infty_\nu(\T, M_d)$. It can be shown that under the pairing:
$$
 [g,f] \bydef \int_\T\operatorname{tr}(f(\zeta)^*g(\zeta))\,dm(\zeta) 
$$ for $f\in L^\infty_\nu(\T, M_d)$ and $g\in L^1_{\nu^*}(\T,M_d)$, the space $L^\infty_\nu(\T, M_d)$ represents the continuous dual space $(L^1_{\nu^*}(\T,M_d))^*$ of $L^1_{\nu^*}(\T,M_d)$ isometrically. This is the usual finite-dimensional vector-valued $L^1$--$L^\infty$ duality. Let
 $$
 \mathcal{H}_-^1(M_d) \bydef \overline{\operatorname{span}} \{ \bar{\zeta}^k T : k \in \N \text{ and } T\in M_d\},
 $$
 where the closure is taken with respect to the $\|\cdot\|_{\nu^*, 1}$ norm. For a subspace $\mathcal{X}$ of the Banach space $L^1_{\nu^*}(\T,M_d)$, we denote its annihilator by $\mathcal{X}^\perp$ which is a closed subspace of $L^\infty_\nu(\T, M_d)$.

 For $\psi = \sum_{j\in \Z} A_j \zeta^j$ in $(\mathcal{H}_-^1(M_d))^\perp$, the space of annihilators of $\mathcal{H}_-^1(M_d)$,
$$
[g, \psi] = \int_\T\operatorname{tr}(\psi(\zeta)^*g(\zeta))\,dm(\zeta) = 0 \quad \text{for all } g\in \mathcal{H}_-^1(M_d).
$$
In particular, 
$$
\int_\T\operatorname{tr}(\bar{\zeta}^j \ \psi(\zeta)^* A)\,dm(\zeta) = 0 \quad \text{for } j \in \N \text{ and } A\in M_d. 
$$
This shows that the
$$
\operatorname{tr}(A_{-j}^* A) = 0 \quad \text{for } j \in \N \text{ and } A\in M_d.
$$
Hence, $A_{-j}=0$ for all $j\in \N$. Thus, the Poisson integral of $\psi$ will produce a function $H^\infty(\D, M_d)$. The above argument can be reversed to show that boundary value of a function in $H^\infty(\D, M_d)$ annihilates $\mathcal{H}_-^1(M_d)$. In conclusion, we have 
\begin{equation}\label{eq:ann}
 ( \mathcal{H}_-^1(M_d))^\perp=H^\infty(\D,M_d).
\end{equation}

The restriction map from the quotient space $(L^1_{\nu^*}(\T,M_d))^*/ (\mathcal{H}_-^1(M_d))^\perp$ to $(\mathcal{H}_-^1(M_d))^*$ given by 
$$
[\cdot, \varphi] + (\mathcal{H}_-^1(M_d))^\perp \mapsto [\cdot, \varphi]|_{\mathcal{H}_-^1(M_d)}, \quad \varphi \in L_\nu^\infty(\T, M_d),
$$
is an isometry; see \cite[Chapter~III, Theorem~10.1]{Conway-book}. Hence, the above $L^1-L^\infty$ duality and the identification \eqref{eq:ann} imply that 
$$
\|\varphi + H^\infty(\D,M_d)\| = \|[\cdot, \varphi]\|_{(\mathcal{H}_-^1(M_d))^*}.
$$
This is a dual version of Nehari-type formula in the present setting which we record below; cf.\cite[Chapter~1, Section~4]{Peller-book}. 
\begin{proposition}\label{prop:dual}
For $\varphi\in L^\infty(\T,M_d)$,
\begin{equation}\label{eq:dual}
 \operatorname{dist}_{L^\infty_\nu}(\varphi, H^\infty(\D, M_d))
 =\sup \{ |[g,\varphi]|: g\in \mathcal{H}_-^1(M_d), \ \|g\|_{\nu^*, 1}\le 1 \}.
\end{equation}
Moreover, the distance is attained.
\end{proposition}

With all this information, let us now return to our interpolation problem with initial nodes $\boldsymbol{\Lambda}=(\lambda_1,\ldots,\lambda_N)$ consisting of distinct points of $\D$ and final nodes $\boldsymbol{W}=(W_1,\ldots,W_N)\in M_d^N$. Consider the finite Blaschke product associated to the initial nodes,
$$
 b_{\boldsymbol{\Lambda}}(z)\bydef \prod_{j=1}^N\frac{z-\lambda_j}{1-\bar{\lambda}_jz},
 $$
 and define
 $$
 \gamma_{\boldsymbol{\Lambda}, \boldsymbol{W}} \bydef \inf \{ \|F\|_{\nu,\infty}: \ F \in H^\infty(\D, M_d) \text{ and } F(\lambda_j)=W_j \text{ for } 1\leq j \leq N \}.
$$
 Consider a function $F_0$ in $H^\infty(\D, M_d)$ that interpolates the data, for example the Lagrange polynomial interpolant 
 $
 \sum_{j=1}^N \ell_j (z)W_j,
 $ where $\ell_j\in \mathbb{C}[z]$ are as in \eqref{LagrangePoly}. Then the collection of all $H^\infty(\D, M_d)$ interpolants can be described via the ideal $b_{\boldsymbol{\Lambda}}H^\infty(\D, M_d)$ as follows:
$$
 \{F\in H^\infty(\D, M_d): F(\lambda_j)=W_j \text{ for } 1\le j \le N\}=F_0+ b_{\boldsymbol{\Lambda}} H^\infty(\D, M_d).
$$
Indeed, for any $H^\infty(\D, M_d)$ interpolant $F$, entrywise division of $F-F_0$ by the finite Blaschke product $b_{\boldsymbol{\Lambda}}$ removes precisely the prescribed zeros and preserves bounded holomorphicity. Since $|b_{\boldsymbol{\Lambda}}|=1$ on $\T$, the classical reduction of interpolation to analytic approximation gives
\begin{align}
 \gamma_{\boldsymbol{\Lambda}, \boldsymbol{W}} &= \inf\{\|F_0 + b_{\boldsymbol{\Lambda}} f\|_{\nu, \infty}: f\in H^\infty(\D, M_d)\}  = \inf\{ \|\bar{b}_{\boldsymbol{\Lambda}} F_0 + f\|_{\nu, \infty} : f\in H^\infty(\D, M_d)\} \notag \\
 &=\operatorname{dist}_{L^\infty_\nu}(\bar{b}_{\boldsymbol{\Lambda}} F_0, H^\infty(\D, M_d)). \label{eq:reduction}
\end{align}
Proposition~\ref{prop:dual} shows that this infimum is attained. Note that if the infumum is achieved by $\psi$ in $H^\infty(\D, M_d)$, then $F_0 - b_{\boldsymbol{\Lambda}} \psi$ is a minimum-norm interpolant.

Motivated by the recent work \cite{Sarkar-KD},  we now utilize the Szeg\"o kernel functions $k_\lambda(z) \bydef (1-\bar{\lambda} z)^{-1}$ on $\D$ for every $\lambda \in \D$ and consider the finite dimensional $M_d$-valued model space
\begin{align}\label{KbLambda}
 \mathcal{K}_{b_{\boldsymbol{\Lambda}}} \bydef  \left\{\sum_{j=1}^Nk_{\lambda_j}C_j:C_j\in M_d\right\}
\end{align}
to obtain a canonical interpolant in $H^\infty(\D, M_d)$. 
Recall the Hilbert-Schmidt inner product on $M_d$ and on $L^2(\T, M_d)$ from Subsection~\ref{sub:basics}.  Let $H^2(\D, M_d)$ be the $M_d$-valued Hardy-Hilbert space on $\D$. For every $g\in H^2(\D, M_d)$,
\begin{align} \label{eq:reproducing}
\langle g, k_\lambda A \rangle_{H^2} = \langle g(\lambda), A\rangle_{M_d} = \operatorname{tr}(A^* g(\lambda)) \quad \text{for } \lambda \in \D , \ A\in M_d. 
\end{align}
Note that $\mathcal{K}_{b_{\boldsymbol{\Lambda}}}$ is the orthogonal complement of $b_{\boldsymbol{\Lambda}} H^\infty(\D, M_d)$ in the Hilbert space $H^2(\T, M_d)$. 
The Gram matrix 
$$
S_{\boldsymbol{\Lambda}} \bydef \left[k_{\lambda_j}(\lambda_i)\right]_{i,j=1}^N
$$
 is positive definite, and hence, is invertible. Then there is a unique $\psi_{\boldsymbol{\Lambda},\boldsymbol{W}}\in  \mathcal{K}_{b_{\boldsymbol{\Lambda}}}$ satisfying the interpolation condition: $\psi_{\boldsymbol{\Lambda},\boldsymbol{W}}(\lambda_j) = W_j$ for $1\leq j \leq N$. More precisely,
\begin{align}\label{UniqueInt}
 \psi_{\boldsymbol{\Lambda},\boldsymbol{W}} =\sum_{j=1}^N k_{\lambda_j}C_j,
 \quad C_j=\sum_{i=1}^N(S_{\boldsymbol{\Lambda}}^{-1})_{ji}W_i,
\end{align}which is a non-zero function provided $\boldsymbol{W}\neq 0$, i.e., at least one of the target matrices $W_j$ is non-zero.

 For $\boldsymbol{W}\ne0$, define the following subspace of $L_{\nu^*}^1(\T, M_d)$:
$$
 \mathcal M_{\boldsymbol{\Lambda},\boldsymbol{W}} \bydef ( \mathcal{K}_{b_{\boldsymbol{\Lambda}}} \ominus \C\psi_{\boldsymbol{\Lambda},\boldsymbol{W}})+ \mathcal{H}_-^1(M_d),
$$
where the orthogonal complement is taken under the Hilbert-Schmidt inner product in $L^2(\T, M_d)$.
The subspace $\mathcal M_{\boldsymbol{\Lambda},\boldsymbol{W}}$ is closed, being a finite-dimensional extension of $\mathcal{H}_-^1(M_d)$. We are now all set to state and prove the Nehari-type distance criterion for our numerical interpolation problem.

\begin{theorem}\label{thm:distance-criterion}
Given an interpolation problem $\boldsymbol{\Lambda} \mapsto \boldsymbol{W}$  with $\boldsymbol{W}\ne0$, let the subspace $\mathcal{K}_{b_{\boldsymbol{\Lambda}}}$ of $H^\infty(\mathbb{D},M_d)$ be as in \eqref{KbLambda} and the unique interpolant $\psi_{\boldsymbol{\Lambda},\boldsymbol{W}}$ in $\mathcal{K}_{b_{\boldsymbol{\Lambda}}}$ be as in \eqref{UniqueInt}. Then 
\begin{align}\label{eq:criterion}
d_{\boldsymbol{\Lambda},\boldsymbol{W}} \bydef \operatorname{dist}_{L_{\nu^*}^1} \left(\frac{\psi_{\boldsymbol{\Lambda},\boldsymbol{W}}}{\|\psi_{\boldsymbol{\Lambda},\boldsymbol{W}}\|_2^2},\mathcal M_{\boldsymbol{\Lambda},\boldsymbol{W}}\right) =\frac1{\gamma_{\boldsymbol{\Lambda},\boldsymbol{W}}}.
\end{align}
Thus, the interpolation problem: $\boldsymbol{\Lambda} \mapsto \boldsymbol{W}$ is solvable by a $\nu$-Schur function exactly when $d_{\boldsymbol{\Lambda},\boldsymbol{W}}\ge1$. Moreover, the problem is extremal precisely when $d_{\boldsymbol{\Lambda},\boldsymbol{W}}=1$ and non-extremal when $d_{\boldsymbol{\Lambda},\boldsymbol{W}}>1$.
\end{theorem}

\begin{proof}
We first show that the closed subspace $\mathcal{E}= \mathcal{K}_{b_{\boldsymbol{\Lambda}}} + \mathcal H_-^1(M_d) \subset L_{\nu^*}^1(\T, M_d)$ is also given by
\begin{equation}\label{eq:test}
 \mathcal{E}= b_{\boldsymbol{\Lambda}}\mathcal H_-^1(M_d).
\end{equation}
Note that multiplication by the Blashcke function $b_{\boldsymbol{\Lambda}}$ on $\mathcal{H}^1_-(M_d)$ or on $L^\infty_\nu(\mathbb{T},M_d)$ is an isometry; thus the space on the right hand side above is a closed subspace.

To establish \eqref{eq:test}, one makes crucial use of the identification \eqref{eq:ann} and the reproducing identity \eqref{eq:reproducing} to conclude that
\begin{equation}\label{eq:test1}
\mathcal {E}^\perp = b_{\boldsymbol{\Lambda}}H^\infty(\D, M_d).
\end{equation}
On the other hand, for every $g\in \mathcal{H}_-^1(M_d)$ and $h\in H^\infty(\D, M_d)$, 
$$
[b_{\boldsymbol{\Lambda}}g, \ h]=[g, \ \bar{b}_{\boldsymbol{\Lambda}}h],
$$ 
and therefore 
\begin{equation}\label{eq:test2}
(b_{\boldsymbol{\Lambda}} \mathcal{H}_-^1(M_d))^\perp=b_{\boldsymbol{\Lambda}} H^\infty(\D, M_d).
\end{equation}
By a standard consequence of the classical Hahn-Banach theorem, the only way two closed subspaces have the same annihilator is that they are the same; thus \eqref{eq:test} is established by \eqref{eq:test1} and \eqref{eq:test2}.

Consider the linear functional on $\mathcal{E}=\mathcal{K}_{b_{\boldsymbol{\Lambda}}} + \mathcal H_-^1(M_d)=b_{\boldsymbol{\Lambda}}\mathcal H_-^1(M_d)$ given by
$$
\ell(g) \bydef  [g,\psi_{\boldsymbol{\Lambda},\boldsymbol{W}}] \quad \text{for all } g \in \mathcal{E}.
$$
Applying \eqref{eq:reduction} for $F_0=\psi_{\boldsymbol{\Lambda},\boldsymbol{W}}$ and then using both \eqref{eq:dual} and \eqref{eq:test} we get
\begin{align} \label{eq:min-norm}
\gamma_{\boldsymbol{\Lambda},\boldsymbol{W}}=\operatorname{dist}_{L^\infty_\nu}(\bar{b}_{\boldsymbol{\Lambda}} \psi_{\boldsymbol{\Lambda},\boldsymbol{W}}, H^\infty(\D, M_d))
&=\sup\{|[b_{\boldsymbol{\Lambda}}g,\psi_{\boldsymbol{\Lambda},\boldsymbol{W}}]|: g\in \mathcal{H}_-^1(M_d), \ \|g\|_{\nu^*, 1}\le1 \} \notag \\
&=\|\ell\|_{\mathcal{E}^*}.
\end{align}
Next, let us note that $\ell$ vanishes on $\mathcal{H}_-^1(M_d)$ by \eqref{eq:ann} and therefore the linear functional $\ell$ on $\mathcal{E}$ boils down to the linear functional $\langle\cdot,\psi_{\boldsymbol{\Lambda},\boldsymbol{W}}\rangle$ restricted to $\mathcal{K}_{b_{\boldsymbol{\Lambda}}}$. Thus,
$$
 \operatorname{ker} \ell=\mathcal M_{\Lambda,W},\quad
 \ell(\psi_{\boldsymbol{\Lambda},\boldsymbol{W}}/\|\psi_{\boldsymbol{\Lambda},\boldsymbol{W}}\|_2^2)=1.
$$
Therefore, a standard application of Hahn-Banach theorem again yields
$$
 \operatorname{dist}_{L_{\nu^*}^1}(\psi_{\boldsymbol{\Lambda},\boldsymbol{W}}/\|\psi_{\boldsymbol{\Lambda},\boldsymbol{W}}\|_2^2, \mathcal{M}_{\boldsymbol{\Lambda},\boldsymbol{W}})=\frac1{\|\ell\|_{\mathcal{E}^*}};
$$
see \cite[Chapter~III, Corollary~6.8]{Conway-book}. 
Using \eqref{eq:min-norm} we finally obtain \eqref{eq:criterion}.
The remaining assertion follows from the definition of $\gamma_{\boldsymbol{\Lambda}, \boldsymbol{W}}$.
\end{proof}


\subsection{Illustrative examples}
We end this section with a couple of solvable interpolation problems -- one is extremal and the other is not. Note that we have seen an example of non-solvable interpolation problem in Example \ref{exp:roadblock}.

\begin{example}\label{Exp-1}
The problem
\begin{equation}\label{eq:data}
 (\lambda_1,\lambda_2)=(0, 4/5)\mapsto
 \left(\begin{bmatrix}0&0\\1&0\end{bmatrix},
 \begin{bmatrix}0&1\\0&0\end{bmatrix}\right)=(W_1,W_2)
\end{equation}
is extremal and hence, there are no choice for the outer functions $\Phi_1,\Phi_2$ as in part (\textbf{NE}) of Theorem~\ref{thm:interpolation}.

The problem is solvable because the rational function
\begin{equation}\label{eq:solution}
 F_0(z)=\frac{1}{5(2-z)^2}
 \begin{bmatrix}
 0&9z\\
 (5-4z)(4-5z)&0
 \end{bmatrix}
\end{equation}is an interpolant. Indeed, it is easy to see that the function $F_0$ is holomorphic in a neighborhood of $\overline{\mathbb{D}}$ (its only pole is at $z=2$), and that it interpolates.  Moreover, from the general fact that
\begin{align} \label{eq:numerical-norm-antidiagonal}
    \nu \left( \begin{bmatrix}
        0 & a \\ b & 0
    \end{bmatrix} \right) = \frac{|a|+|b|}{2}, \ \ a, b \in \mathbb{C}
\end{align} it follows that for $\zeta \in \mathbb{T}$
\begin{align*}
    \nu(F_0(\zeta)) = \frac{9 + |5-4\zeta| |4-5\zeta|}{10 |2-\zeta|^2} = 1.
\end{align*}
Hence, by \eqref{RadLim-Iso}, $\|F_0\|_{\nu, \infty} =1$.

To see that this problem is extremal, let $F=(F_{jk})_{1\leq j, k \leq 2}$ be a $\nu$-Schur function which interpolates the data. We show that $\|F\|_{\nu, \infty}=1$. 

For $x=\frac{1}{\sqrt{2}}(1, i)^{\textup{T}}$ and $y=\frac{1}{\sqrt{2}}(1, -i)^{\textup{T}}$ in $\C^2$, we find that 
$$
i(f_{12}(z)-f_{21}(z)) = \langle f(z)x, x \rangle - \langle f(z)y, y \rangle \quad \text{for } z\in \D.
$$
Then the holomorphic function $g(z)=(f_{12}(z)-f_{21}(z))/2$ on $\D$ satisfies:
\begin{align}\label{eq:g}
  |g(z)| \leq \nu(f(z)) \leq 1 \quad \text{for } z\in \D.  
\end{align}
Moreover, $g(0)=-1/2$ and $g(4/5)=1/2$. So, the holomorphic function 
$$
h(z)= \frac{g(z)+1/2}{1+ g(z)/2}
$$
maps $\D$ to $\D$ with $h(0)=0$ and $h(4/5)=4/5$. Therefore, by the Schwarz lemma we must have $h(z)=z$ for every $z\in \D$ and hence, 
$$
g(z)= \frac{z-1/2}{1-z/2} \quad \text{for } z\in \D.
$$
Thus, $\|g\|_\infty =1$. Therefore, \eqref{eq:g} implies that $\|F\|_{\nu, \infty} =1$.
\end{example}

We now turn to an example of a non-extremal interpolation problem. 
\begin{example} \label{Exp-2}
 The interpolation problem
$$
\bm\Lambda=\left(-\frac12, 0, \frac12 \right) \mapsto \left( \begin{bmatrix} 0 & \frac38\\ 0 & 0
 \end{bmatrix}, \begin{bmatrix} 0&0\\ 0&0 \end{bmatrix}, \begin{bmatrix} 0 & \frac38\\ 0 & 0\end{bmatrix}\right) = \textbf{W}
 $$
is {\em non-extremal}. In fact, the argument below shows that
$$
\gamma_{\bm\Lambda, \boldsymbol{W}}
=\inf\left\{
      \|F\|_{\nu, \infty}: F\in H^\infty(\D, M_d),\quad
      F(\lambda_j)=W_j \ \text{ for } 1\leq j\leq 3
    \right\}=\frac34.
$$ 
The rational function 
$$ f(z)= \begin{bmatrix} 0&\frac32z^2\\\ 0&0 \end{bmatrix}
$$
clearly, interpolates the data and 
 $\|f\|_{\nu, \infty} = \sup_{z\in \D}\frac34|z|^2 = \frac34$.

We can factor $f$ as $f(z)=P(z)G(z)Q(z)$, where
$$
P(z) = \begin{bmatrix}
    \frac{\sqrt{7}}{2} & 0\\ 0 & \frac{1}{2}
\end{bmatrix}, \quad G(z) = \begin{bmatrix}
    0 & \frac67 z^2 \\0 & 0
\end{bmatrix},\quad \text{and} \quad Q(z)=\begin{bmatrix}
 \frac{1}{2}  & 0\\ 0 & \frac{\sqrt{7}}{2} 
\end{bmatrix}.
$$
Clearly, $G$ is a rational Schur function with $\|G\|_\infty = \frac{6}{7}$ and the polynomials $P$ and $Q$ are invertible satisfying $PP^* + Q^*Q = 2 I_2$.

We now show that if $F$ is any interpolant, then $\|F\|_{\nu,\infty} \geq \frac{3}{4}$. Towards this, consider the unit vector $x=\frac1{\sqrt2}(1, 1)^T \in \C^2$ and define
$$
g(z)= \langle F(z)x,x\rangle \quad \text{for } z\in \D.
$$
Then $g$ is holomorphic and $\|g\|_\infty \leq \|F\|_{\nu,\infty}$. The interpolation conditions imply
$$g(0)=0,\quad g\left(\frac12\right) =g\left(-\frac12\right)=\frac3{16}.$$
The holomorphic function $$h(z) = \frac{g(z)+g(-z)}2$$
therefore satisfies $\|h\|_\infty \leq \|F\|_{\nu,\infty}$, $h(0)=0$ and $h'(0)=0$.
By Schwarz lemma 
we get $|h(z)|\leq \|F\|_{\nu,\infty} |z|^2$ for each $z\in \D$.
Putting $z=1/2$ in this inequality we get
$$\frac3{16} =\left|h\left(\frac12\right)\right| \leq\frac{\|F\|_{\nu,\infty}}{4}
$$which gives the required inequality $\|F\|_{\nu,\infty}\geq3/4$.
\end{example}

\section{The $\nu$-inner functions}\label{sec:inner}
We devote this section to the study of $\nu$-inner functions, a special class of $\nu$-Schur functions. Recall the notation $\overline{\Omega}_\nu$ for the closed numerical radius unit ball from \eqref{Nu-Ball}. The $\nu$-inner functions maps $\T$ into a specific subset of the topological boundary of $\overline{\Omega}_\nu$, which we define next.
\begin{definition}\label{def:inner}
A numerical contraction $T$ in $M_d$ is said to be {\it extremal numerical contraction} if the only possible way an operator of the form
$$ \begin{bmatrix} T & B \\ 0 & C \end{bmatrix}: \begin{bmatrix}
      \mathbb{C}^d \\ \mathcal K
 \end{bmatrix}\to \begin{bmatrix}
      \mathbb{C}^d \\ \mathcal K
 \end{bmatrix}
$$is a numerical contraction is $B=0$. We shall denote by $\partial_e$ the collection of all extremal numerical contractions.

We denote the collection of extremal numerical contractions by $\partial_e$.
\end{definition}
In \cite{ABDS}, the authors provided various algebraic and operator theoretic characterizations of extremal numerical contractions. The next lemma, proved in \cite[Proposition~2.1]{ABDS}, provides one such description via the set of common positive defect:
\begin{definition}\label{D:Defect}
For a numerical contraction $T$, define its set of {\it common positive defects} by
\begin{equation*} 
 \mathcal D(T) \bydef
 \{R \in M_d^{\rm sa} : 0 \preceq R \preceq I_d - \operatorname{Re}(\zeta T) \text{ for every } \zeta \in \mathbb{T}\}.
\end{equation*}
\end{definition}

Clearly, $0 \in \mathcal D(T)$. Moreover, integrating $R \preceq I_d - \operatorname{Re}(\zeta T)$ over $\mathbb{T}$ we have $0 \preceq R \preceq I_d$. Thus, $\mathcal D(T)$ is a nonempty, closed, convex, and compact subset of $M_d^{\textup{sa}}$. 

In particular, the lemma below shows that $\partial_e$ is a subset of the topological boundary of $\overline{\Omega}_\nu$, which constitutes numerical contractions of $\nu$-norm one.
\begin{lemma}\label{lem:defect}
A numerical contraction $T$ is in $\partial_e$ if and only if $\mathcal D(T) = \{0\}$. In particular, every extremal numerical contraction satisfies $\nu(T) = 1$.
\end{lemma}
We need the following simple lemma.
\begin{lemma} \label{lem:complex-extreme}
Let $T$ be a numerical contraction. For every $R\in\mathcal D(T)$,
\begin{equation}\label{eq:complex-disc}
 T+zR\in \Omega_{\nu} \quad \text{for all } z \in \overline{\D}.
\end{equation}
\end{lemma}
\begin{proof}
For any \( R \in \mathcal{D}(T)\) and \(z \in \overline{\D}\), we have 
$$
 I_d - \operatorname{Re}(\zeta(T+zR))
 =I_d - \operatorname{Re}(\zeta T)-\operatorname{Re}(\zeta z)R
 \succeq I_d - \operatorname{Re}(\zeta T)-R\succeq0.
$$This completes the proof.
\end{proof}
For the next characterization of extremal numerical contractions, let us denote
$$
(M_d^+)_1 \bydef \{Z \in M_d^{\textup{sa}}: \ 0 \preceq Z \preceq I_d\}
$$and define $\varphi: (M_d^+)_1 \to \mathbb{R}$ by
\begin{align}\label{KiFunctionEta}
\varphi(Z) \bydef \|Z\|_{\textup{HS}}^2 - \operatorname{tr} Z,  \ \ Z\in (M_d^+)_1.
\end{align}
\begin{lemma}\label{lem:extremal-Borel}
Let $T$ be a numerical contraction and $\varphi$ be the map as in \eqref{KiFunctionEta}. Then $\varphi|_{\mathcal{D}(T)}$ attains its minmum exactly at one point, denoted by $ \Phi(T)$. The map $T \mapsto \Phi(T)$ from $\overline{\Omega}_\nu$ to $(M_d^+)_1$, is Borel measurable.
Moreover, $\Phi(T) \neq 0 $ if and only if $T \notin \partial_e$. 
\end{lemma}

\begin{proof}
Let $T\in \overline{\Omega}_\nu$. Then $\varphi|_{\mathcal{D}(T)}$ has a minimum on the non-empty compact set $\mathcal{D}(T)$. If possible assume that it has two distinct minimizers $R_1, R_2$ in $\mathcal{D}(T)$. By the parallelogram identity in $M_d$ endowed with the tracial inner product, we have 
\begin{align*}
\varphi\left(\frac{R_1 + R_2}{2}\right)= \frac{\|R_1\|_{\textup{HS}}^2 + \|R_2\|_{\textup{HS}}^2}{2} -\frac14 \|R_1 - R_2\|_{\textup{HS}}^2 -  \frac12 \operatorname{tr} (R_1 + R_2 ) &< \frac{ \varphi(R_1) + \varphi(R_2)}{2} \\ 
& = \min \varphi|_{\mathcal{D}(T)}
\end{align*}
The above strict inequality is absurd as $(R_1 + R_2)/2$ belongs to the convex set $\mathcal{D}(T)$.
Thus, $\Phi(T)$ is the unique minimizer of $ \varphi$ on $\mathcal{D}(T)$.
So, Lemma~\ref{lem:defect} tell us that 
$T\notin \partial_e$ if and only if $\Phi(T) \neq 0$.

Also, $ \varphi$ satisfies the hypothesis of Lemma~\ref{lem:Borel} with the compact subset
$$
\mathcal{G} =\{(T, R) \in \overline{\Omega}_\nu \times (M_d^+)_1: \ R \in \mathcal{D}(T) \} 
$$
of $\overline{\Omega}_\nu \times (M_d^+)_1$. Hence, we conclude that $\Phi$ is Borel measurable.
\end{proof}

We are now ready to define the inner functions in the present setting.
\begin{definition}\label{def:inner_func}
A function $f \in H^\infty(\mathbb{D}, M_d)$ is said to be {\it $\nu$-inner} if the associated radial limit function $f^\flat$ takes values in $\partial_e$ almost everywhere on $\T$, the set of extremal numerical contractions.
\end{definition}
It is clear that every unitary matrix is an extremal numerical contraction, and thus, the standard inner matrix-valued functions are $\nu$-inner functions. For a more non trivial example, note that the function 
$$
f(z) = \begin{bmatrix} 0 & 2z \\ 0 & 0 \end{bmatrix}
$$ is $\nu$-inner because whenever $z\in\mathbb{T}$, the matrix above is an extremal numerical contraction.

\begin{remark}
In the classical setting, a function in $H^\infty(\mathbb{D}, M_d)$ with sup-norm bounded by one is called \emph{inner} if almost everywhere on $\mathbb{T}$, it takes values in the set of unitary matrices, which constitute the Shilov boundary of the operator-norm unit ball of $M_d$.

The Shilov boundary of the unit ball $\overline{\Omega}_\nu$ is completely characterized in the authors' companion manuscript \cite{ABDS}: while extremal numerical contractions form part of this Shilov boundary, they do not exhaust it. In view of this, a reader might regard our definition of an inner function as somewhat restrictive, since $\nu$-inner functions take values in a proper subset of the Shilov boundary. Nevertheless, we show that every non-extremally solvable $\nu$-interpolation problem admits a $\nu$-inner solution, which in turn yields the density of $\nu$-inner functions in the $\nu$-Schur class.
\end{remark}

\begin{lemma}\label{lem:measurable-defect}
Let \(\psi:\T\to \overline{\Omega}_{\nu}\) be a Borel measurable function. If \(\psi(\zeta)\) is not
extremal numerical contraction on a set of positive measure, there is a Borel measurable function \(h:\T\to M_d\) such that for every $\zeta \in \T$,
\begin{equation}\label{eq:measurable-defect}
 0\preceq h(\zeta)\preceq I_d,\quad
 h(\zeta)\in\mathcal D(T(\zeta)),\quad \text{and} \quad 
 \|h\|_{\textup{HS}, 2}>0.
\end{equation}
\end{lemma}

\begin{proof}
Define $h \bydef \Phi \circ \psi$ with the Borel map $\Phi$ as in Lemma~\ref{lem:extremal-Borel}. Then $h$ is Borel as it is a composition of two Borel maps. Also, $h(\zeta)\in\mathcal{D}(\psi(\zeta))\subset (M_d^+)_1$ for every $\zeta$ in $\T$. Also, by Lemma~\ref{lem:extremal-Borel},
\[
\mathcal F \bydef \{\zeta\in\T:\ \psi(\zeta)\text{ is not extremal}\}=\{\zeta \in \T:\ h(\zeta)\ne0\},
\]
which is a Borel set. By our hypothesis, $\mathcal F$ has $m(\mathcal F)>0$. The function $\zeta\mapsto\|h(\zeta)\|_{\textup{HS}}^2$ is Borel, nonnegative, and strictly positive on $\mathcal F$. Therefore,$\|h\|_{\textup{HS}, 2}^2=\int_\T \|h\|_{\textup{HS}}^2\,dm\ge\int_{\mathcal F}\|h\|_{\textup{HS}}^2\,dm > 0$, because a nonnegative function that is strictly positive on a set of positive measure has positive integral.
\end{proof}
In classical scalar and matrix-valued Nevanlinna--Pick interpolation, rational inner interpolants---specifically, finite Blaschke and Blaschke--Potapov products---are constructively produced via the Schur algorithm~\cite{Ball-Gohberg-Rodman1990, DGK-1, DGK-2}. This recursive scheme depends crucially on the transitive action of the automorphism group of the matrix operator-norm unit ball. As noted in the Introduction, however, the automorphism group of the numerical radius unit ball fails to act transitively: as shown in the authors' forthcoming paper~\cite{ABDS}, every such automorphism is linear, leaving the orbit of the zero matrix a singleton and obstructing any Schur-algorithmic approach. A second standard avenue for constructing inner interpolants is the state-space realization method, which has proved effective on the unit disk, the polydisk~\cite{AM-book, Knese}, and quotients of the polydisk~\cite{BK-JFA, DKS}. Because an analogous realization theory for $\nu$-Schur functions is not yet available---a development we leave to a forthcoming paper---this operator-theoretic route is likewise unavailable here. Nevertheless, relying directly on function-theoretic machinery, we prove below that every non-extremal numerical radius interpolation problem admits a $\nu$-inner solution.

\begin{theorem}\label{thm:strict}
Let $\bf\Lambda\mapsto \boldsymbol{W}$ be a non-extremal interpolation problem. Let $\mathcal{S}_{\bf\Lambda,\boldsymbol{W}}$ denote the set of $\nu$-Schur interpolants of this data. Then the set of $\nu$-inner functions intersected with  $\mathcal{S}_{\bf\Lambda,\boldsymbol{W}}$ contains a dense $G_\delta$ subset of $\mathcal{S}_{\bf\Lambda,\boldsymbol{W}}$ in the compact-open topology.
\end{theorem}
\begin{proof}
By assumption on the data, there exists $f_0$ such that $f_0 \in \mathcal{S}_{\bf\Lambda, \boldsymbol{W}}$ and 
\begin{align}\label{delta}
\|f_0\|_{\nu,\infty} \le 1-\delta
\end{align}for some $0<\delta <1$. An application of Montel's theorem shows that the set $\mathcal{S}_{\bf\Lambda,\boldsymbol{W}}$ is a nonempty compact metric space in the compact-open topology of the spaces of $M_d$-valued holomorphic functions in $\D$. Consider the boundary-value map:
\[
 J:\mathcal{S}_{\bf\Lambda,\boldsymbol{W}} \longrightarrow L^2(\mathbb{T},M_d) \quad\mbox{given by}\quad J(f)=f^\flat.
\]
For each fixed \(0\le r<1\), the map
\[
 J_r:\mathcal{S}_{\bf\Lambda,\boldsymbol{W}} \longrightarrow L^2(\mathbb{T},M_d) \quad\mbox{given by} \quad J_r(f)(\zeta)=f_r(\zeta)=f(r\zeta)
\]
 for $\zeta \in \mathbb{T}$, is continuous as the locally uniform convergence is uniform on the circle
\(|z|=r\). Moreover, for $f, g \in \mathcal{S}_{\bf\Lambda,\boldsymbol{W}}$,
\[
 \|J_r(f)-J_r(g)\|_{\textup{HS}, 2}
 \le\sqrt d\,\sup_{|z|=r}\|f(z)-g(z)\|.
\]
This is obtained using the fact that $\|T\|_{\textup{HS}} \leq \sqrt{d} \|T\|$ for $T\in M_d$.
Further, we have
$$
f_r = J_r(f)\xrightarrow{r\to 1^-} f^\flat= J(f) \text{ in } L^2(\mathbb{T}, M_d),
$$
for every fixed $f$; cf. \cite[Theorem~4.7.B]{RR-book}. Choosing a sequence $\{r_n\}_n$ of positive real numbers increasing to $1$, we have $J_{r_n}$ converging to $J$ pointwise. Thus, by application of Baire-Osgood theorem for functions between two complete metric spaces, we infer that the set $\mathcal{G}$ consisting of points of continuity of the map $J$ forms a dense $G_\delta$ subset of $\mathcal{S}_{\bf\Lambda,\boldsymbol{W}}$; cf. \cite[Theorem~11.20]{Carothers}

It will be enough if we prove that every $f\in\mathcal G$ is $\nu$-inner. Let us fix $f\in\mathcal G$. Suppose on the contrary that $f$ is not $\nu$-inner, i.e., by Definition \ref{def:inner_func}, the set
\begin{align}\label{BadSet}
\{\zeta\in \mathbb{T}: f^\flat(\zeta)=\lim_{r\to 1-}f(r\zeta)\mbox{ exists but not extremal}\}
\end{align}
is a set of positive measure.

Our goal is to construct a sequence $F_k$ of functions in $\mathcal{S}_{\bm\Lambda, \textbf{W}}$ such that $\nu(F_k)<1$ for each $k$ and $F_k\to f$ uniformly over compact subsets of $\mathbb{D}$, yet $J(F_k)$ does not converge to $J(f)$, contradicting the continuity of $J$ at $f$. The construction of the sequence will put us through some rigorous analysis.

For each $0<r<1$, $f_r$ is continuous on $\mathbb{T}$ and hence in particular, is Borel measurable. The set $E$ where the radial limit $\lim_{r\to 1-}f(r\zeta)$ exists is therefore Borel measurable. We assign $f^\flat(\zeta)=0$ whenever $\zeta\notin E$. The resulting function, which we continue to denote by $f^\flat$, is therefore a Borel measurable function with values in $\overline{\Omega}_\nu$. By Lemma~\ref{lem:measurable-defect} applied to $f^\flat$ we obtain a Borel measurable function $h$ such that 
\begin{equation}\label{eq:h-properties}
0<\|h\|_{\textup{HS}, 2}<\infty,\quad 0\preceq h(\zeta)\preceq I_d \quad\mbox{and}\quad h(\zeta)\preceq I_d-\operatorname{Re}(\omega f^\flat(\zeta))
\end{equation}for every $\zeta,\omega \in \mathbb{T}$.

For \(0<r<1\), we apply the  Poisson integral \eqref{Poisson} on both sides of \eqref{eq:h-properties} to get the following for every $\zeta,\omega\in\mathbb{T}$: 
\begin{align*}
 0\preceq h(r\zeta) = \int_\mathbb{T} \frac{1-r^2}{|\eta - r\zeta|^2} h(\eta) \ dm(\eta) 
 &
 \preceq \int_\mathbb{T} \frac{1-r^2}{|\eta - r\zeta|^2} \big(I_d - \operatorname{Re}(\omega f^\flat(\zeta)\big) \ dm(\eta)  
\\
 &\preceq
 I_d-\operatorname{Re}(\omega f(r\zeta)).
\end{align*}
Thus, $h_r(\zeta) \bydef h(r\zeta) \in \mathcal{D}(f(r\zeta))$ for every $0<r<1$, and so, Lemma~\ref{lem:complex-extreme} implies 
\begin{equation}\label{eq:poisson-disc}
 \nu(f(r\zeta)+z h_r(\zeta))\le1
 \quad \text{for every }\zeta\in\mathbb{T},\ z\in \overline{\mathbb{D}}.
\end{equation}
 We next convert these measurable defect directions into analytic perturbations that preserve all the prescribed values.

With the Lagrange polynomials
$$
 \ell_j(z) \bydef \prod_{\substack{1\le k\le N\\ k\ne j}}
                \frac{z-\lambda_k}{\lambda_j-\lambda_k},
 \quad 1\le j\le N
$$
let us define the matrix-valued polynomials
\[
 c^{(r)}(z) \bydef \sum_{j=1}^N\ell_j(z)\bigl(W_j-f(r\lambda_j)\bigr)
\]for every $0<r<1$. Continuity of $f$ at the finitely many nodes \( \lambda_1, \ldots, \lambda_N\) and
$$
 \|c^{(r)}\|_\infty \le\sum_{j=1}^N
    \sup_{|z|\le1}|\ell_j(z)|\,\|W_j - f(r\lambda_j)\|
$$
shows that \(\|c^{(r)}\|_\infty\to0\) as \(r \to 1^-\). Therefore, with $0<\delta<1$ as in \eqref{delta}, we have
for every $k\geq1,$ there exists an $k/(k+1)< r_k <1$ such that 
\begin{equation}\label{eq:correction-small}
\|c^{(r_k)}\|_\infty < \delta/4(k+1)
\end{equation}
Using the Lagrange identity $\ell_j(\lambda_i)=\delta_{i,j}$ we see that for every $k\geq1$ and $j=1,2,\dots, N$,
\begin{equation}\label{eq:corrected-values}
 f(r_k\lambda_j)+c^{(r_k)}(\lambda_j)=W_j.
\end{equation}

Note that with the Borel measurable function $h$ as in \eqref{eq:h-properties}, $h_{r_k}$ is a continuous $M_d$-valued function on $\mathbb{T}$ and so, $h_{r_k}$ can be approximated uniformly in $\mathbb{T}$ by matrix-valued Laurent polynomials. Thus, we find an $M_d$-valued Laurent polynomial
$$
 p_k(\zeta)=\sum_{j=-L_k}^{L_k}A_{k,j}\zeta^j
$$
such that for every $k\geq1$,
\begin{equation}\label{eq:Laurent-small}
 \sup_{\zeta\in\mathbb{T}}\|p_k(\zeta)-h_{r_k}(\zeta)\|
       \le\frac{\delta}{4(k+1)}.
\end{equation} 
 The norm inequality 
$$
 \|p_k-h\|_{\textup{HS},2} \leq \|p_k - h_{r_k}\|_{\textup{HS},2} + \|h_{r_k} - h\|_{\textup{HS},2} \leq \sqrt{d}\ \|p_k - h_{r_k}\|_{\infty,\mathbb{T}}+ \|h_{r_k} - h\|_{\textup{HS},2},
$$
the equation \eqref{eq:Laurent-small} and the fact that \(\|h -  h_r\|_{\textup{HS},2} \to 0\) as \(r\to 1^-\) together imply 
\begin{equation}\label{eq:p-L2}
 \|p_k-h\|_{\textup{HS},2} \to 0\quad \mbox{as} \quad k\to\infty.
\end{equation}

Let the initial data set be $\bm\Lambda=(\lambda_j)_{j=1}^N$, and consider the finite Blaschke product
$$
 b_{\bf\Lambda}(z) = \prod_{j=1}^N\frac{z-\lambda_j}{1-\overline{\lambda_j}z}.
$$
Since $z^n\to 0$ uniformly over compact subsets of $\D$, we see that for each $k\geq1$, we can find a sufficiently large integer \(M_k>L_k\) such that
\begin{equation}\label{eq:oscillation-small}
 \sup_{|z|\le \frac{k}{k+1}}
   \|b_{\bf\Lambda}(z)z^{M_k}p_k(z)\| \leq 
   \sup_{|z|\le \frac{k}{k+1}}
   \|b_{\bf\Lambda}(z)z^{L_k}p_k(z)\| \cdot \sup_{|z|\le \frac{k}{k+1}} \|z^{M_k - L_k}\| \le\frac1k.
\end{equation}

With the interpolant $f_0$ as in \eqref{delta}, let us define the sequence of holomorphic functions
\begin{equation}\label{eq:Fk}
 F_k(z)\bydef \frac{k}{k+1} \bigl(f(r_kz)+c^{(r_k)}(z)+b_{\bf\Lambda}(z)z^{M_k}p_k(z)\bigr) + \frac{1}{k+1} f_0(z).
\end{equation}
What is clear is that each $F_k$ is in $H^\infty(\mathbb{D}, M_d)$. The following steps will lead us to the desired contradiction.
\smallskip

\noindent \textbf{Claim 1:} For each $k\geq1$, $F_k$ solves the interpolation problem $\bm\Lambda\mapsto \textbf{W}$ and $\nu(F_k)<1$. In particular, $F_k \in \mathcal{S}_{\bm\Lambda, \textbf{W}}$ for each $k\geq1$.
\smallskip

\noindent \textbf{Claim 2:} The functions $F_k$ converges uniformly over compact subsets of $\mathbb{D}$ to the function $f$, where $f$ is the fixed function in the set $\mathcal{G}$ of continuity of the function $J$ with the property \eqref{BadSet}.
\smallskip

\noindent \textbf{Claim 3:} But $J(F_k)$ does not converge to $J(f)$ as $k\to\infty$, contradicting the continuity of $J$ at $f$.
\smallskip

\noindent {\sf{Proof of Claim 1:}} The fact that each $F_k$ interpolates the data follows from equations
\eqref{eq:corrected-values} and the facts that $f_0$ interpolates and \(b_{\bf\Lambda}(\lambda_j)=0\). We note that for almost every $\zeta \in \T$, the radial limit function $F_k^\flat$ of $F_k$ can be written as
\begin{align*}
 F_k^\flat(\zeta) &= \frac{k}{k+1} \bigg(f(r_k\zeta)+c^{(r_k)}(\zeta)+b_{\bf\Lambda}(\zeta)\zeta^{M_k}p_k(\zeta)\bigg) +\frac{1}{k+1}f_0^\flat(\zeta)\\
 &
=\frac{k}{k+1}
       \bigg(f(r_k\zeta)+b_{\bf\Lambda}(\zeta)\zeta^{M_k}h_{r_k}(\zeta)\bigg)
       +\frac{1}{k+1} f_0^\flat(\zeta)+ \\
       &\quad \quad \quad 
       \frac{k}{k+1}
       \bigg(c^{(r_k)}(\zeta)
          +b_{\bf\Lambda}(\zeta)\zeta^{M_k}(p_k(\zeta)-h_{r_k}(\zeta))\bigg).
\end{align*}
This together with the inequalities in \eqref{delta},
\eqref{eq:poisson-disc},
\eqref{eq:correction-small}, and \eqref{eq:Laurent-small} show
\begin{align}\label{eq:Fk-bound}
 \nu(F_k^\flat(\zeta))\le \frac{k}{k+1}+\frac{1}{k+1}(1-\delta)
       +\frac{k}{k+1}\frac{\delta}{2(k+1)}\le1-\frac{\delta}{2(k+1)}.
\end{align}
From the general observation made in \eqref{RadLim-Iso}, it follows that 
$
\|F_k\|_{\nu,\infty} = \|F_k^\flat\|_{\nu,\infty} < 1 .
$
\smallskip

\noindent {\sf{Proof of Claim 2:}}
 We have observed the following uniform convergences in the compact-open topology:  (i) \(f_r \to f \) as $r \to 1^-$; (ii) \(c^{(r_k)}\to0\) as $k\to\infty$, by \eqref{eq:correction-small}; and (iii) $b_{\bf\Lambda}(z)z^{M_k}p_k(z)\to 0$ as $k\to\infty$, by
\eqref{eq:oscillation-small}. Finally, since \(f,f_0\) both are bounded, equation \eqref{eq:Fk} gives $F_k\to f$
in the compact-open topology.
\smallskip

\noindent {\sf{Proof of Claim 3:}}
Let us denote for each $k\geq1$
\[
 a_k(z) \bydef \frac{k}{k+1}\bigg(f(r_kz)+c^{(r_k)}(z)\bigg)
                +\frac{1}{k+1}f_0(z).
\]Then for each $k\geq1$, we have
\begin{align}\label{AlternateExpression}
F_k(z) = a_k(z) + \frac{k}{k+1}b_{\bm\Lambda}(z)z^{M_k}p_k(z).
\end{align} In the following computation, we obtain the first inequality by adding and subtracting $\frac{1}{k+1} f$ to $a_k-f$ and then use the triangle inequality:
\begin{align}\label{eq:ak-L2}
  \|a_k-f \|_{\textup{HS}, 2} \notag
  & \le   \frac{k}{k+1} \|f_{r_k} - f + c^{(r_k)}\|_{\textup{HS}, 2} + \frac{1}{k+1} \| f_0 -f\|_{\textup{HS}, 2} \\
  & \leq \frac{k}{k+1} \|f_{r_k} - f\|_{\textup{HS}, 2} + \frac{k}{k+1}\|c^{(r_k)}\|_{\textup{HS}, 2} + \frac{1}{k+1} \| f_0 -f \|_{\textup{HS}, 2} ,
\end{align}
which converges to zero as $k$ tends to $\infty$. The reverse triangle inequality and the expression \eqref{AlternateExpression} of $F_k$ is used to obtain
\[
 \|F_k - f \|_{\textup{HS}, 2}
 \ge \frac{k}{k+1}\|p_k\|_{\textup{HS}, 2}-\|a_k -f \|_{\textup{HS}, 2}.
\]
Together with \eqref{eq:p-L2}, \eqref{eq:ak-L2}, and
\eqref{eq:h-properties}, this yields
\begin{equation}\label{eq:contradiction}
 \liminf_{k\to\infty}\|J(F_k)-J(f)\|_{\textup{HS}, 2}\ge\|h\|_{\textup{HS}, 2}>0.
\end{equation}
This is a contradiction because $f$ is a point of continuity of $J$ and $F_k$ converges to $f$ in the domain $\mathcal{S}_{\bf\Lambda, \boldsymbol{W}}$ of $J$.

Thus, $f^\flat$ must take value in $\partial_e$, the set of extremal numerical contractions almost everywhere on $\mathbb{T}$. Since $f\in\mathcal G$ was arbitrary, $\mathcal G$ is contained in the set of $\nu$-inner functions intersected with $\mathcal{S}_{\bf\Lambda, \boldsymbol{W}}$.
The density and $G_\delta$ property of $\mathcal G$ completes the proof.
\end{proof}

Finally, we establish the density of $\nu$-inner functions in the $\nu$-Schur class under the compact-open topology.

\begin{theorem}\label{thm:main}
For every $F \in H^\infty(\mathbb{D}, M_d)$ with $\|F\|_{\nu,\infty} \le 1$, there exists a sequence of $\nu$-inner functions $\{F_k\}$ converging to $F$ in the compact open topology. 
\end{theorem}

\begin{proof}
Let $F$ be as in the statement. We choose two sequences:
$$
 \lambda_j=\frac1{j+2} \quad(j\ge1),\, \text{ and }\, 
 t_k=\frac{k}{k+1}\quad(k\ge1).
$$
The nodes $\{\lambda_j : j\in \mathbb{N}\}$ are distinct and accumulate at \(0\in\mathbb{D}\).
For each \(N\geq1\), let us consider the interpolation problem
\[
\bm\Lambda_N = (\lambda_1, \lambda_2,\dots,\lambda_N)\longmapsto (t_N F(\lambda_1),t_N F(\lambda_2),\dots, t_N F(\lambda_N)) = \textbf{W}_N.
\]
By choice, these problems are solvable by strict $\nu$-Schur functions, viz., \(t_NF\). By Theorem~\ref{thm:strict}, there exist $\nu$-inner functions $F_N$ that solves the problem $\bm\Lambda_N\mapsto \textbf{W}_N$ for each $N\geq1$, i.e.,
\begin{equation}\label{eq:scaled-values}
 F_N(\lambda_j)=t_NF(\lambda_j) \quad \mbox{for each } j=1,2,\dots, N.
\end{equation}
Since $\|F_N\|_{\infty} \le 2\|F_N\|_{\nu,\infty} < 2$ for each $N\geq1$, Montel's theorem gives a convergent subsequence $\{F_{N_k}\}_{k\geq1}$ converging to $F_0$ (say) in the compact open topology of $H^\infty(\mathbb{D}, M_d)$.

Let us now note that given any fixed $j \in \mathbb{N}$, we have a large enough $k$ depending on $j$ so that \(N_k \ge j\). Evaluating
\eqref{eq:scaled-values} at \(\lambda_j\) and using \(t_{N_k}\to1\) as $k\to\infty$ gives 
$$
F_0(\lambda_j)= \lim_{k\to\infty} F_{N_k}(\lambda_j) = \lim_{k\to\infty}
t_{N_k} F(\lambda_j) = F(\lambda_j).
$$
Since the sequence accumulates at a domain point, viz., the origin, $F_0$ must be identically equal to $F$ by the identity theorem for scalar holomorphic functions applied to each entry of the matrix functions under considerations.

The argument above actually shows more. Indeed, by applying Montel's theorem to every subsequence of $\{F_k\}$, we see that every subsequence of $\{F_k\}$ has a convergent subsequence whose limit coincides with $F$. Thus the whole sequence $\{F_k\}_k$ must converge to $F$ in the compact-open topology.
\end{proof}

\vspace{2mm}

\noindent\textbf{Acknowledgement:} The first named author is supported by Prime Minister's Research Fellowship, ID: 1303115. The research of the second author is supported by ANRF-NPDF (PDF/2025/002836). The third and fourth named author's research is supported by the Anusandhan National Research Foundation (ANRF), Government of India, under the Mathematical Research Impact-Centric Support (MATRICS) scheme Grant No. ANRF/ARGM/2025/002455/MTR and ANRF/ARGM/2025/000573/MTR, respectively.

The fourth author is grateful to Joseph A.~Ball for valuable discussions that led to the necessary condition in Theorem~\ref{thm:Nec-cond} in 2017.
\vspace{2mm}

\noindent\textbf{AI disclosure:}
Although this project was set in motion long before AI tools became prevalent, the authors acknowledge the use of various large language models (LLMs) as research aids during the preparation of this article. Specifically, the strategy for obtaining a Borel-measurable selection of self-adjoint operators for the Ando completion of a $\nu$-Schur function via Lemma~\ref{lem:Borel}, which is subsequently employed in Theorem~\ref{thm:gen-structure} and Lemma~\ref{lem:measurable-defect}, was developed with the assistance of LLMs. Examples~\ref{exp:roadblock}, \ref{Exp-1}, and~\ref{Exp-2} were constructed with LLM assistance. The manuscript was written entirely by the authors, who have independently verified every argument and take full responsibility for the mathematical correctness of the paper.


\begin{thebibliography}{99}
\bibitem{Abrahamse1979}
Abrahamse, M. B. {\em A Pick interpolation theorem for finitely connected domains}. Michigan Math. J. 26 (1979), no. 2, 195--203.

\bibitem{Abrahamse-Douglas1976}
Abrahamse, M. B.; Douglas, R. G. {\em A class of subnormal operators related to multiply connected domains}. Adv. Math. 19 (1976), no. 1, 106--148.

\bibitem{Agler-OpThAdvAppl}
Agler, Jim. {\em On the representation of certain holomorphic functions defined on a polydisc}. Topics in operator theory: Ernst D. Hellinger memorial volume, 47–66. Oper. Theory Adv. Appl., 48. Birkhäuser Verlag, Basel, 1990.
ISBN:3-7643-2532-1.

\bibitem{AM-Crelle} Agler, Jim; McCarthy, John E.
{\em Nevanlinna-Pick interpolation on the bidisk}.
J. Reine Angew. Math. 506 (1999), 191–204.


\bibitem{AM-book} Agler, Jim; McCarthy, John E. {\em Pick Interpolation and Hilbert Function Spaces}. American Mathematical Society, Providence, 2002.

\bibitem{ABDS} Ahmed, H.; Bhowmik, M.; Das, B. Krishna; Sau, H. {\em The Shilov boundary and automorphisms of the numerical radius ball}. Preprint in final stage of preparation.

\bibitem{Ando} Ando, T. {\em Structure of operators with numerical radius one}.
Acta Sci. Math. (Szeged) 34 (1973), 11–15.

\bibitem{Ball-Bolotnikov} Ball, Joseph A.; Bolotnikov, Vladimir. {\em Nevanlinna-Pick interpolation for Schur-Agler class functions on domains with matrix polynomial defining function in $\C^n$}. 
New York J. Math. 11 (2005), 247–290.



\bibitem{Ball-Clancey1992}
Ball, Joseph A.; Clancey, Kevin F. {\em Interpolation by $H^\infty$ functions on a Riemann surface}. J. Funct. Anal. 107 (1992), no. 2, 295--313.

\bibitem{Ball-Gohberg-Rodman1990}
Ball, Joseph A.; Gohberg, Israel; Rodman, Leiba. {\em Interpolation of Rational Matrix Functions}. Operator Theory: Advances and Applications, 45. Birkh\"auser Verlag, Basel, 1990.

\bibitem{Ball-Helton1983}
Ball, Joseph A.; Helton, J. William. {\em A Beurling-Lax theorem for the Lie group $U(m,n)$ which contains most classical interpolation theory}. J. Operator Theory 9 (1983), no. 1, 107--142.

\bibitem{Ball-Trent1998}
Ball, Joseph A.; Trent, Thomas T. {\em The Nevanlinna-Pick problem for the Bergman space and related Hilbert function spaces}. J. Funct. Anal. 157 (1998), no. 1, 1--61.


\bibitem{Barclay-JLMS} Barclay, Steven. {\em Continuity of the spectral factorization mapping}.
J. London Math. Soc. (2) 70 (2004), no. 3, 763–779.

\bibitem{Bebiano} Bebiano, Natalia; Spitkovsky, Ilya M. {\em Numerical ranges of Toeplitz operators with matrix symbols}.
Linear Algebra Appl. 436 (2012), no. 6, 1721–1726.

\bibitem{Bercovici} Bercovici, Hari; Foias, Ciprian; Tannenbaum, Allen. {\em A spectral commutant lifting theorem}. Trans. Amer. Math. Soc. 325 (1991), no. 2, 741–763.

\bibitem{Bhatia} Bhatia, Rajendra; Jain, Tanvi. {\em The numerical radius and positivity of block matrices}.
Linear Algebra Appl. 656 (2023), 463–482.

\bibitem{BK-JFA}Bhowmik, Mainak; Kumar, Poornendu. {\em Function theory on quotient domains related to the polydisc}.
J. Funct. Anal. 289 (2025), no. 6, Paper No. 110978, 29 pp.

\bibitem{Caratheodory} Caratheodory, C. {\em Theory of functions of a complex variable}. Vol. 2.
Translated by F. Steinhardt. Chelsea Publishing Co., New York, 1954. 220 pp.

\bibitem{Carothers} Carothers, N. L. {\em Real analysis}. Cambridge University Press, Cambridge, 2000, xiv+401 pp.
ISBN: 0-521-49756-6.

\bibitem{Choi} Choi, Man-Duen; Li, Chi-Kwong. {\em Numerical ranges and dilations}. Linear and Multilinear Algebra 47 (2000), no. 1, 35–48.

\bibitem{Conway-book} Conway, John B. {\em A course in functional analysis}.
Grad. Texts in Math., 96.
Springer-Verlag, New York, 1990, xvi+399 pp.
ISBN: 0-387-97245-5.

\bibitem{DKS} Das, B. Krishna; Kumar, Poornendu; Sau, Haripada. {\em Distinguished varieties and the Nevanlinna-Pick problem on the symmetrized bidisk}.
Math. Z. 310 (2025), no. 3, Paper No. 40, 29 pp.


\bibitem{DGK-1} Delsarte, Ph.; Genin, Y.; Kamp, Y.  {\em Schur parametrization of positive definite block-Toeplitz systems}.
SIAM J. Appl. Math. 36 (1979), no. 1, 34–46.

\bibitem{DGK-2}Delsarte, Ph.; Genin, Y.; Kamp, Y. {\em The Nevanlinna-Pick problem for matrix-valued functions}.
SIAM J. Appl. Math. 36 (1979), no. 1, 47–61.

\bibitem{DM} Dritschel, Michael A.; McCullough, Scott. {\em Test functions, kernels, realizations and interpolation}.
Operator theory, structured matrices, and dilations, 153–179.
Theta Ser. Adv. Math., 7. Theta, Bucharest, 2007.
ISBN: 978-973-87899-0-6.

\bibitem{DW} Dritschel, Michael A.; Woerdeman, Hugo J. {\em Model theory and linear extreme points in the numerical radius unit ball}. Mem. Amer. Math. Soc. 129 (1997), no. 615, viii+62 pp.

\bibitem{FFGK}Foias, C.; Frazho, A. E.; Gohberg, I.; Kaashoek, M. A. {\em Metric constrained interpolation, commutant lifting and systems}.
Oper. Theory Adv. Appl., 100.
Birkhäuser Verlag, Basel, 1998, xii+587 pp.
ISBN: 3-7643-5889-0.

\bibitem{Gustafson-Rao1997}
Gustafson, Karl E.; Rao, D. K. M. {\em Numerical Range: The Field of Values of Linear Operators and Matrices}.
Universitext, Springer-Verlag, New York, 1997.

\bibitem{Helson-book} Helson, Henry. {\em Lectures on invariant subspaces}. Academic Press, New York-London, 1964. xi+130 pp.

\bibitem{HL}Helson, Henry; Lowdenslager, David.
{\em Prediction theory and Fourier series in several variables}.
Acta Math. 99 (1958), 165–202.

\bibitem{Sarkar-KD} K. D., Deepak; Sarkar, J. {\em Commutant lifting, interpolation, and perturbations on the polydisc}.
To appear in the American Journal of Mathematics. arXiv:2301.10020v3.


\bibitem{Knese} Knese, Greg. {\em Schur-Agler class rational inner functions on the tridisk}.
Proc. Amer. Math. Soc. 139 (2011), no. 11, 4063–4072.


\bibitem{Nevanlinna} Nevanlinna, R. {\em \"Uber beschränkte Funktionen, die in gegebenen Punkten vorgeschrieben Werte annehmen}. Ann. Acad. Sci. Fenn. 13 No. 1 (1919), 1–72.

\bibitem{Pick} Pick, Georg. {\em Über die Beschränkungen analytischer Funktionen, welche durch vorgegebene Funktionswerte bewirkt werden}. Math. Ann. 77 (1915), no. 1, 7–23.

\bibitem{Wiener-Masani-I} Masani, P.; Wiener, N. {\em The prediction theory of multivariate stochastic processes. I. The regularity condition}. Acta Math. 98 (1957), 111–150.

\bibitem{Peller-book} Peller, Vladimir V. {\em Hankel operators and their applications}.
Springer Monogr. Math. Springer-Verlag, New York, 2003, xvi+784 pp. ISBN: 0-387-95548-8.

\bibitem{Rosenblatt}Rosenblatt, Murray. {\em A multi-dimensional prediction problem}. Ark. Mat. 3 (1958), 407–424.

\bibitem{RR-book} Rosenblum, Marvin; Rovnyak, James. {\em Hardy classes and operator theory}.
Oxford Math. Monogr. Oxford University Press, New York, 1985.

\bibitem{Sarason1967}
Sarason, Donald. {\em Generalized interpolation in $H^\infty$}. Trans. Amer. Math. Soc. 127 (1967), no. 2, 179--203.


\end{thebibliography}
\end{document}